\pdfoutput=1
\RequirePackage{ifpdf}
\ifpdf % We~are running pdfTeX in pdf mode
\documentclass[pdftex]{sigma}
\else
\documentclass{sigma}
\fi

\numberwithin{equation}{section}

\newtheorem{Theorem}{Theorem}[section]
\newtheorem*{Theorem*}{Theorem}
\newtheorem{Corollary}[Theorem]{Corollary}
\newtheorem{Lemma}[Theorem]{Lemma}
\newtheorem{Proposition}[Theorem]{Proposition}
\newtheorem{Conjecture}[Theorem]{Conjecture}

\theoremstyle{definition}
\newtheorem{Definition}[Theorem]{Definition}

\newtheorem{Remark}[Theorem]{Remark}

\newcommand{\eps}{\varepsilon}

\newcommand{\FF}{{\mathbb{F}}}

\newcommand{\R}{{\mathbb R}}
\newcommand{\CC}{{\mathbb C}}

\newcommand{\cP}{{\mathcal P}}

\newcommand{\lstar}{{\raise-0.15ex\hbox{$\scriptstyle \ast$}}}

\newcommand{\ind}{{\bf{1}}}
\newcommand{\cB}{\mathcal{B}}

\newcommand{\cF}{\mathcal{F}}

\newcommand{\cS}{\mathcal{S}}
\newcommand{\cT}{\mathcal{T}}

\newcommand{\Tr}{\textup{Tr}}
\def\eqd{\stackrel{d}{=}}

\newcommand{\sgn}{\operatorname{sgn}}
\newcommand{\Pf}{\operatorname{Pf}}
\newcommand{\Hf}{\operatorname{Hf}}
\newcommand{\ww}{\boldsymbol\omega}
\newcommand{\nn}{\boldsymbol\eta}
\newcommand{\cA}{\mathcal{A}}
\newcommand{\cG}{\mathcal{G}}
\newcommand{\cD}{\mathcal{D}}
\newcommand{\dd}{\Theta}

\newcommand{\T}{\dag}

\newcommand{\lst}[1]{[\![#1 ]\!]}
\newcommand{\nint}[2]{\lfloor #1 \rfloor_{#2}}

\newcommand{\nfr}[2]{\left\{ #1 \right\}_{#2}}

\newcommand{\mbf}[1]{\mathbf{#1}}

\newcommand{\wt}[1]{\widetilde{#1}}

\newcommand{\HH}{\mathtt{H}_{n,\beta}}
\newcommand{\WW}{\mathtt{W}_{n,m, \beta}}
\newcommand{\SQW}{\mathtt{SqW}_\beta}

\begin{document}

\allowdisplaybreaks

\newcommand{\arXivNumber}{2412.04579}

\renewcommand{\PaperNumber}{049}

\FirstPageHeading

\ShortArticleName{Solvable Families of Random Block Tridiagonal Matrices}

\ArticleName{Solvable Families of Random Block Tridiagonal\\ Matrices}

\Author{Brian RIDER~$^{\rm a}$ and Benedek VALK\'O~$^{\rm b}$}

\AuthorNameForHeading{B.~Rider and B.~Valk\'o}

\Address{$^{\rm a)}$~Department of Mathematics, Temple University, Philadelphia, PA, USA}
\EmailD{\mail{brian.rider@temple.edu}}

\Address{$^{\rm b)}$~Department of Mathematics, University of Wisconsin -- Madison, Madison, WI, USA}
\EmailD{\mail{valko@math.wisc.edu}}

\ArticleDates{Received November 05, 2025, in final form April 23, 2026; Published online May 14, 2026}

\Abstract{We introduce two families of random tridiagonal block matrices for which the joint eigenvalue distributions can be computed explicitly. These distributions are novel within random matrix theory, and exhibit interactions among eigenvalue coordinates beyond the typical mean-field log-gas type. Leveraging the matrix models, we go on to describe the point process limits at the edges of the spectrum in two ways: through certain random differential operators, and also in terms of coupled systems of diffusions. Along the way we establish several algebraic identities involving sums of Vandermonde determinant products.}


\Classification{60B20; 15A15}

\section{Introduction}\label{sec1}

Trotter observed that if one applies the Householder tridiagonalization process to a GOE or GUE random matrix then the resulting real symmetric tridiagonal matrix will have independent entries (up to symmetry) with normal and chi distributions \cite{Trotter}. In~\cite{DE}, Dumitriu and Edelman presented a far reaching generalization of this result. They show that, for any $\beta > 0$, the~${ n \times n}$ random Jacobi matrix with independent $N\bigl(0,\frac{2}{\beta}\bigr)$ random variables along the diagonal, and independent \smash{$ \frac{1}{\sqrt{\beta}} \chi_{\beta(n-1)}, \frac{1}{\sqrt{\beta}} \chi_{\beta(n-2)}, \dots, \frac{1}{\beta} \chi_\beta$} random variables along the off-diagonals, has joint eigenvalue density proportional to
\begin{equation}
\label{eig_DE}
\left|\Delta(\lambda)\right|^\beta {\rm e}^{-\frac{\beta}{4} \sum_{j=1}^n \lambda_j^2}.
\end{equation}
Here $\Delta(\lambda)$ denotes the usual Vandermonde determinant of the eigenvalues.
This includes Trotter's result for GOE or GUE upon setting $\beta=1$ or $2$.

The Dumitriu--Edelman model for the Gaussian, or ``Hermite'', beta ensemble, along with their Laguerre counterparts, initiated an immense amount of activity in the study of the scaling limits of beta ensembles. See, for instance,
\cite{ES, KillipNenciu,KS,KRV,RR, RRV,BVBV,BVBV_sbo}.
Motivated both by the original construction of \cite{DE} along with its ensuing impact, here we establish two families of similarly solvable block-tridiagonal matrix models.

Let $\HH (r,s)$ denote the distribution of the $rn \times rn$ symmetric or Hermitian block tridiagonal matrix with $r \times r$ diagonal blocks
distributed as independent copies of
G(O/U)E, and descending upper diagonal blocks distributed as independent copies of
the (lower triangular) positive square root of a real/complex Wishart with parameters $(r, (r+s)(n-i))$. Here $i$ is the index of the offdiagonal block entry, and
$\beta=1$ and 2 corresponds to the real and complex case, respectively.
As in the $r=1$ case, the diagonal and offdiagonal variables are also independent of each other. A more detailed description of these ensembles is provided in Section \ref{subs:matrix_distr}.

Note of course that the Wishart distribution is the multivariate analog of the $\chi^2$ distribution, and that $\HH(1,s)$ is just the original Dumitriu--Edelman model, after a reparameterization.
Further, when $s=0$, our model may be arrived at by a suitable block tridiagonalization procedure of the corresponding $rn \times rn$ G(O/U)E, {\`a} la Trotter. This has already been noticed in~\cite{Spike2} in the context of eigenvalue spiking. Finding a suitable general beta version of the spiked Tracy--Widom laws introduced in that paper was another motivation for our work. Modifying the parameter of the offdiagonal block square root Wishart distributions in a similar manner to that of the offdiagonal $\chi$ distributions in the Dumitriu--Edelman model presented a natural starting point.

Our main result is the following.

\begin{Theorem}\label{thm:main}
For $\beta =1$ and $2$, the symmetrized joint eigenvalue density of $\HH(r,s)$ can be computed explicitly in the following cases:
 \begin{align}
 \frac{1}{Z_{n, \beta, r, 2}} |\Delta({\lambda})|^{\beta}
 \left( \sum_{(\mathcal{A}_1,\dots,\mathcal{A}_r)\in \cP_{r,n}} \prod_{j=1}^r \Delta(\cA_j)^2 \right) {\rm e}^{- \frac{\beta}{4}\sum_{i=1}^{rn} \lambda_i^2} \qquad
 \mbox{for}\ r \ge 2, \ \beta s=2, \label{density1}
 \end{align}
and
 \begin{align}
 \label{density2}
 \frac{2^n}{Z_{n, \beta, 2, \beta s}} \Delta({\lambda})^{\beta+\frac{\beta s}{2}} \left|\Pf \left(\frac{{\bf{1}}_{i \neq j}}{\lambda_i -\lambda_j} \right)\right|^{\frac{\beta s}{2}} {\rm e}^{- \frac{\beta}{4}\sum_{i=1}^{2n}\lambda_i^2} \qquad \mbox{for} \ r = 2, \ \beta s = 2,4.
 \end{align}
It further holds that
\begin{align*}
Z_{n, \beta, r, \beta s} ={}&(n r)! (2\pi)^{\frac{nr}{2}}
\left(\frac{\beta}{2}\right)^{a_{n,\beta,r,s}} \Gamma\left(\frac{\beta}{2}\right)^{-nr} \prod_{k=1}^{nr} \Gamma\left(\frac{\beta}{2}\left(k+s \left\lceil\frac{k}{r}\right\rceil\right)\right)\\
& \times \begin{cases}
 1, &\beta s=2,\\
 \left(\dfrac{\beta}{12}\right)^n, &\beta s=4,
\end{cases}
\end{align*}
with $a_{n,\beta,r,s}=
-\frac{\beta}{4} n r (n (r+s)+s)+\bigl(\frac{\beta}{4}-\frac{1}{2}\bigr){nr}$
for all $n$, $\beta = 1$ and $2$, and combinations of $r$ and $s$ in \eqref{density1} and \eqref{density2}.
\end{Theorem}
Here
for $r\ge 2$ and $n\ge 1$, $\cP_{r,n}$ denotes the set of size $r$ equipartitions of $\lst{rn} := \{ 1,2, \dots, rn\}$.
That is, $\{\cA_1, \dots,\cA_r\}\in \cP_{r,n}$ if $|\cA_i|=n$ for all $i$ and the $\cA_i$ form a partition of $\lst{rn}$. With that, for any
 $\cA \subset \lst{rn}$, we write $\Delta(\cA)$ as shorthand for the Vandermonde determinant in the $|\cA|$ ordered eigenvalue variables with indices drawn from $\cA$ (suppressing the explicit dependence on~$\lambda_i$, $i \in \cA$). Finally,
$\Pf(M)$ denotes the Pfaffian of $M$.

In both \eqref{density1} and \eqref{density2}, we see novel types of interactions among the points beyond the usual $|\Delta({\lambda})|$
to some power.
The formulas for the overlapping $r=2$, $\beta s = 2$ cases are shown to agree by a Pfaffian/Vandermonde identity, see Lemma~\ref{lem:det4_identities} below.
This is one of several identities involving sums of powers of Vandermonde determinants that we prove in Section \ref{sec:det_identities}.
We also note that \eqref{density1} is consistent with \eqref{eig_DE} upon taking $r=1$, as then the sum over equipartitions reduces to $\Delta(\lambda)^2 = \Delta(\lambda)^{\beta s}$. Of course, unlike the $r=1$ case, our results are restricted to a special set of values of $s$. Further, the parameter $s$ does not appear to have a natural interpretation, as an inverse temperature say, in the
resulting density formulas.

One might anticipate that the form of the $r=2$ family should generalize to all even integer $\beta s$. However, computer assisted calculations for small $n$ values indicate that the Pfaffian structure in~\eqref{density2} breaks down for $\beta s=6$. Understanding what happens for larger block size $r$ beyond~${\beta s=2}$ also remains open. Our difficulty in extending exact formulas to either parameter regime is
tied to our approach to proving Theorem~\ref{thm:main}. This rests on computing the absolute $\beta s$-moment of a certain structured determinant over the Haar distributed Orthogonal or Unitary group (in dimension $rn$). We do this by expansion and re-summation, the underlying complexity of which grows in both $r$ and $\beta s$. In another direction,
our block model could certainly be constructed using quaternion ingredients, leading to $\HH(r,s)$ with $\beta=4$.
The non-commutativity of the quaternion variables poses additional technical challenges in extending Theorem~\ref{thm:main} to that setting, though we expect these are not insurmountable.

Next, a natural question is whether densities of the form \eqref{density1} or \eqref{density2} appear ``in the wild". In fact,
the $r=2$ family bears close resemblance to what is known as the
Moore--Read, or Pfaffian, state for the fractional quantum Hall effect, see \cite{MR_1991}. In that theory the points lie in the complex plane,
so \eqref{density2} might be viewed as a one-dimensional caricature of these states in the same way that the Gaussian (and other) beta ensembles are one-dimensional caricatures of a~true coulomb gas. Additionally, the sum-over-equipartitions interaction term in \eqref{density1} is by itself reminiscent of the structure found for the densities of eigenvalue powers of certain determinantal point processes, see, for instance, \cite{Rains}.

The eigenvalues of random block matrices have of course been studied in a number of capacities, most notably perhaps as structured band matrices connected to the Anderson or Wegner orbital models, see, for example,~\cite{SchSch}
and the references therein. Motivated by the theory of matrix orthogonal polynomials,
\cite{Dette1} and \cite{Dette2} introduce families of ``block beta'' Hermite, Laguerre and Jacobi ensembles built out of Gaussian and/or $\chi$ variables, and study their limiting density of states.
The large deviations of related ensembles have been considered in~\cite{Rouault1,Rouault2}.
Our work though is the first to provide a systematic approach to finding solvable block models.

We close the introduction with descriptions of: (i) the soft edge asymptotics for $\HH(r,s)$, and (ii), how the results stated through that point, including the associated asymptotics, extend to a family of block Wishart (or Laguerre) ensembles. After this, Section \ref{sec2} lays out some basic facts on the spectral theory of block tridiagonal matrices along with the detailed definitions of our various matrix models. Section \ref{sec3} provides an overview of the eigenvalue density derivations, identifying a certain moment calculation as fundamental (see Theorem~\ref{thm:moment}). That calculation is spread over Sections \ref{sec:2moment} and \ref{sec:4moment}, for moments $\beta s =2$ and $\beta s = 4$ respectively. Section \ref{sec:det_identities} establishes a number of identities (and presents a conjecture in a related spirit) involving
sums of Vandermonde determinant powers required in the preceding. Finally,
 Section \ref{sec:asymptotics} is devoted to asymptotics.

\subsection[Soft edge asymptotics of HH(r,s)]{Soft edge asymptotics of $\boldsymbol{\HH(r,s)}$}

While it does not appear possible to compute correlations directly from the formulas \eqref{density1} or~\eqref{density2}, the random operator approach is available. In the block setting, this was developed by Bloemendal and Vir\'ag for the soft edge in \cite{Spike2}, and their approach applies to our case for any values of $r$ and $s$. It even applies in the $\beta=4$ case where we do not have statements about joint eigenvalue densities.

Introduce the $\beta =1,2,$ or $4$ matrix Brownian motion $B_x$ in dimension $r$:
 the independent, stationary increment process for which $B_y- B_x \sim B_{y-x}$ is distributed as $\sqrt{y-x}$ times a copy of $r \times r$ G(O/U/S)E.
 Next, for $\gamma > 0$, bring in the differential operator acting on $r$-dimensional vector valued functions on $\R_{+}$,
\begin{equation}\label{eq:H_op}
 \mathcal{H}_{\beta, \gamma} = - \frac{{\rm d}^2}{{\rm d}x^2} + rx + \sqrt{\frac{2}{\gamma}} B'_x.
\end{equation}
When $\gamma=1$, this is the multivariate stochastic Airy operator of~\cite{Spike2}. In particular, with a~Dirichlet boundary condition at the origin, the spectrum of $-\mathcal{H}_{\beta} = -\mathcal{H}_{\beta, 1}$ is given by the $\operatorname{Airy}_\beta$ process, the edge scaling limit of the Gaussian beta ensemble. The largest value of this process (which is minus the ground state eigenvalue of $\mathcal{H}_{\beta}$), has classical Tracy--Widom distribution~$TW_\beta$ with~$\beta =1,2, 4$.

\begin{Theorem}
\label{thm:limit_op}
For any $r$, $s$ and $\beta=1,2,4$,
let $\mathbf{T}_n \sim \HH(r,s)$. Denote by \smash{$\lambda_0^{(n)} < \lambda_1^{(n)} < \cdots $} the eigenvalues
of the renormalized
\begin{equation*}
\mathbf{H}_n = \gamma^{-1/2} (rn)^{1/6} \bigl(2 \sqrt{(r+s)n} {I}_{rn} - \mathbf{T}_n \bigr),
\end{equation*}
and by $\Lambda_0 < \Lambda_1 < \cdots$ the Dirichlet eigenvalues of
$ \mathcal{H}_{\beta, \gamma}$ with the choice $\gamma = \frac{r+s}{r}$ .
Then the point process \smash{$\big\{ \lambda_0^{(n)} ,\lambda_1^{(n)} , \dots\big\}$} converges in law to $\{\Lambda_0, \Lambda_1, \dots \}
$ as $n\to \infty$.
\end{Theorem}

Implicit in the above is that the edge of the spectrum of $\HH(r,s)$ occurs around $2 \sqrt{(r+s) n}$. In fact, and not surprisingly, the normalized counting measure of eigenvalues of $\HH(r,s)$ converges to a scaled semi-circle law. A statement (and proof) of this fact is given below as Theorem~\ref{thm:eigcount}.

Theorem~\ref{thm:limit_op} itself follows from the main result of \cite{Spike2}. We do though sketch an overview of the ideas in Section~\ref{sec:asymptotics}. Similarly, \cite[Theorem~1.5]{Spike2} provides a second description of the limiting point process $\{ \Lambda_i \}_{i \ge 0}$ via matrix oscillation theory. Applying the same here yields the following.

\begin{Corollary}
\label{cor:osc}
Define the measure $\mathbb{P}$ on non-intersecting paths $\mbf{p}=(p_1, \dots, p_r)\colon[0,\infty) \mapsto ( -\infty, \infty]$ induced by the stochastic differential equation system
\begin{equation}
\label{mult_sde}
 {\rm d}p_i = \frac{2}{\sqrt{\beta \gamma}} {\rm d}b_i + \bigg(\lambda + rx - p_i^2 + \sum_{j \neq i} \frac{2}{p_i - p_j} \bigg){\rm d}x,\qquad 1\le i \le r,
\end{equation}
starting from $(p_1(0), \dots , p_r(0)) = \{\infty\}^r$
and entering $\{ p_1 < \cdots < p_r\}$ at $x>0$. Here $(b_1, \dots, b_k)$ is a standard real $r$-dimensional Brownian motion; $p_1$ can hit $-\infty$ in finite time,
whereupon it is placed at $+\infty$ and the re-indexed process starts afresh. Then with $\Lambda_0< \Lambda_1< \cdots $ defined as in Theorem~{\rm\ref{thm:limit_op}}, and again $\gamma = \frac{r+s}{r}$, it holds that
$ P( \Lambda_k \le \lambda ) = \mathbb{P}$ $(x\mapsto \mbf{p}(x)$ explodes at most ${k}$ times$)$
for all $k \ge 0$.
\end{Corollary}

The above corollary immediately implies that, whenever $\beta \gamma$ equals a classical value, i.e., $1$, $2,$ or $4$, we can deduce that the limiting edge point process
corresponds to that of the
G(O/U/S)E. In particular, in this case $\Lambda_0$ will have $TW_{\beta \gamma}$ distribution. This again is one of the primary takeaways of \cite{Spike2}.
Due the fact that the eigenvalue law at finite $n$ is independent of any
chosen
block tridiagonalization (independent of the choice of $r$), it follows that, again when the diffusion parameter is classical, the explosion times of \eqref{mult_sde} are equal in distribution for all $r\ge 1$. No direct proof of this striking fact is known.

Specifying to the cases for which we have novel explicit joint eigenvalue densities, this implies the following.

\begin{Corollary}
\label{cor:betalimit}
Consider the random point process defined by the $r=2$, $\beta s = 2$ joint density~\eqref{density1} in Theorem~{\rm\ref{thm:main}}.
When $\beta=1$, the
appropriately rescaled
point process converges in law to the $\operatorname{Airy}_2$ point process. In the case of $r=2$ and $\beta s= 4$, the appropriately scaled process determined by \eqref{density2} in Theorem~{\rm\ref{thm:main}} converges in law to the $\operatorname{Airy}_4$
point process when $\beta=2$. In particular, in these cases the largest eigenvalues $($after rescaling$)$ converge to the classical $TW_2$ and $TW_4$ distributions, respectively.
\end{Corollary}

Conjecturing that the $r$-fold diffusion characterization of Corollary \ref{cor:osc} provides the description of the $\operatorname{Airy}_{\beta \gamma}$ process for any $\beta \gamma>0$ we arrive to the following.

\begin{Conjecture}
\label{con:betalimit}
More generally, the point process scaling limit of \eqref{density1} is distributed as $\operatorname{Airy}_{\beta+2/r}$ for all $r \ge 2$ and $\beta =1$ or $2$. In the case of \eqref{density2} with $\beta s = 4$ and $\beta=1$, the point process scaling limit is $\operatorname{Airy}_{3}$.
\end{Conjecture}

\subsection{Block Laguerre ensembles}

In \cite{DE}, the authors also produce $\beta$ generalizations of the classical Laguerre (Wishart) ensemble, showing that there is an $n\times n$ tridiagonal matrix model built out of independent $\chi$ variables for which the
eigenvalue density is proportional to
\begin{equation}
\label{eig_DE1}
\left|\Delta(\lambda)\right|^\beta
\prod_{i=1}^n \lambda_i^{\frac{\beta}{2}(m-n+1) -1} {\rm e}^{-\frac{\beta}{2} \sum_{i=1}^n \lambda_i} \mathbf{1}_{\R_+^n}.
\end{equation}
When $\beta =1$ or $2$, this coincides with that of the law of a sample covariance matrix for $m\ge n$ independent real or complex normal
samples in dimension $n$.
Along with $\beta$ now taking any positive value, the model behind \eqref{eig_DE1} allows $m$ to be generalized to any real number greater than $n-1$.

We define the distribution $\mathtt{W}_{n,m, \beta}(r, s)$ on nonnegative definite block tridiagonals as follows. Let $\mathbf{L}_n$
be an $rn \times rn$ block bidiagonal matrix with independent $r\times r$ diagonal and upper offdiagonal blocks denoted by
$\{\mbf{D}_i\}_{i=1,n}$ and $\{\mbf{O}_i\}_{i=1, n-1}$, that are lower and upper triangular matrices, respectively.
 Distribute these according to
square-root Wishart matrices with parameters
$(r, (r+s)(m+1 -i))$ and $(r, (r+s)(n-i))$, respectively. Then $\mathtt{W}_{n, m, \beta}(r, s)$ is the
law of~$\mbf{L}_n \mbf{L}_n^\dagger$. Full details are provided in Definition \ref{def:BlockW}.

Again, when $s=0$ this model has been considered previously in \cite{Spike2} and \cite{RR} in connection to eigenvalue spiking.
In that case the underlying random matrix $\mbf{L}_n$ arises from an explicit block bi-diagonalization of an $rn \times rm$ matrix of independent Gaussians.

Effectively the same considerations behind Theorem~\ref{thm:main} imply the following.

\begin{Theorem}\label{thm:main_W}
The joint eigenvalue density of $\mathtt{W}_{n, m, \beta}(r, s)$ for $\beta=1$ or $2$ has the form
\eqref{density1} for general $r\ge 2$ and $\beta s=2$ and \eqref{density2} for $r =2$ and $\beta s =2$ or $4$ with an explicitly computable normalizing constant,
the only change being that the Gaussian weight \smash{$ {\rm e}^{-\frac{\beta}{4} \sum_{i=1}^{rn} \lambda_i^2}$} is replaced by~\smash{$\prod_{i=1}^{rn} \lambda_i^{\beta( (r+s)(m-n)+1)/2-1} {\rm e}^{-\beta \lambda_i/2}$}, restricted to $\R_{+}^{rn}$.
\end{Theorem}

In terms of asymptotics, the normalized eigenvalue counting measure will converge to a~Ma\-rchenko--Pastur distribution, under suitable conditions on $n$ and $m$. See again Theorem~\ref{thm:eigcount} below.

Toward accessing something more refined,
we focus on the choice $m = n +a $ for fixed $a > -1/(r+s)$ as $n \rightarrow \infty$ and look at the scaling limit of the smallest eigenvalues, which end up being in the vicinity of the origin. This is the random matrix hard edge, and introduces novel limiting phenomena beyond what we have for $\mathtt{H}_{n, \beta}(r, s)$. Note that, similar to Theorem~\ref{thm:limit_op}, it may be proved that the suitably centered and scaled largest eigenvalues under $\mathtt{W}_{n, m, \beta}(r, s)$ will converge to those of $\mathcal{H}_{\beta, \gamma}$, for an appropriate $\gamma$. The same is true for the smallest eigenvalues when $\liminf_{n\to \infty} m/n>1$.

For the hard edge, the characterizing limit operator is now of Sturm--Liouville type:
again acting on $r$-dimensional vector valued functions,
%$f \in L^2(\R_+, \mathbb{F}^r)$,
\begin{equation}
\label{matrixgenerator}
 \mathcal{G}_{\beta, \gamma} = - {\rm e}^{rx} {\bf{Z}}_x \frac{\rm d}{{\rm d}x} {\mbf{Z}_x^{-1} } \frac{\rm d}{{\rm d}x}.
\end{equation}
Here $x \mapsto {{\mbf{Z}}_x} $ is a symmetrized version of drifted Brownian on the general real or complex linear group dimension $r$, the parameters $\gamma$ and $a$ coefficients of the defining stochastic differential equation (see \eqref{WandA} below).
Similar to $\mathcal{H}_{\beta, \gamma}$, the operator $\mathcal{G}_{\beta, \gamma}$ for $\gamma =1$
has previously been shown to characterize multi-spiked hard edge laws \cite{RR2} for
$\beta =1,2,4$. For $\gamma=1$ and $r=1$, this is the stochastic Bessel operator introduced by Ram\'{\i}rez and Rider in
\cite{RR}.
In analogy to Theorem~\ref{thm:limit_op} and Corollary \ref{cor:osc}, we have the following.

\begin{Theorem}
 \label{thm:limit_op1} For $\mbf{W}_n \sim \mathtt{W}_{ n, n+a, \beta}(r, s)$, denote by
 \smash{$0 < {\lambda}_0^{(n)} < {\lambda}_1^{(n)} < \cdots $} the point process of eigenvalues of \smash{$ \frac{rn}{\gamma} \mbf{W}_n$}.
As $n \rightarrow \infty$ this converges in law to the point process
\smash{$0 < \hat{\Lambda}_0< \hat{\Lambda}_1 <\cdots $}
of
Dirichlet eigenvalues of $ \mathcal{G}_{\beta, \gamma}$
with $\gamma = \frac{r+s}{r}$.
%Here again $\gamma = \frac{r+s}{r}$.
\end{Theorem}

The dependence on the many underlying parameters is made more explicit in the Riccati picture.

\begin{Corollary}
\label{cor:osc1}
Let $\mathbb{P}$ be the measure on non-intersecting paths ${\mathbf{q}}\colon [0, \infty) \mapsto [-\infty, \infty]^r$
defined~by
\begin{equation*}
\label{rrq}
 {\rm d} q_{i} = \frac{2}{\sqrt{\beta \gamma}} q_{i} {\rm d} b_i + \left( \left(\frac{a}{\gamma} + \frac{2}{\beta \gamma}\right) q_{i} - q_{i}^2 - \lambda {\rm e}^{-r x} + q_{i} \sum_{j \neq i} \frac{ q_{i} + q_{j}}{ q_{i}- q_{j} } \right) {\rm d} x,
\end{equation*}
started from $\{ \infty\}^r$ with the same ordering and re-indexing conventions upon possible passages to $-\infty$ described in
Corollary {\rm\ref{cor:osc}}. With $\gamma =\frac{r+s}{r}$ and $0 < \hat{\Lambda}_0< \hat{\Lambda}_1 <\cdots $ defined in Theorem~{\rm\ref{thm:limit_op1}},
it holds
%\label{HardEdge_zeros}
$ P \bigl(\hat{\Lambda}_k > \lambda\bigr) = \mathbb{P}$ $( x \mapsto \mbf{q}(x)$ vanishes at most $ k $ times$)$
for any given $k = 0,1,\dots$.
\end{Corollary}

Again, whenever $\beta \gamma = 1, 2$ or $4$ we conclude that the point process scaling limit of the smallest eigenvalues of ${\mathtt{W}}_{n, n+a,\beta} (r, s)$ is the classical hard edge, or Bessel, point process. More generally,
we conjecture that these limits are given by the general $\beta \gamma$ hard edge process
defined in \cite{RR}. In particular, versions of Corollary \ref{cor:betalimit} and Conjecture \ref{con:betalimit} are readily formulated. We record these at the end of Section \ref{sec:asymptotics}.

 Having dealt with the soft and hard edge scaling limit of our models, it is natural to ask if the same can be done in the bulk case. The analogous results to \cite{Spike2} and \cite{RR2} for the bulk have not though yet been developed.
Another possible future direction is to extend our results to circular ensembles using the results of \cite{KillipNenciu} as a starting point.

\section{Preliminaries}\label{sec2}

We start by outlining some basic facts on the spectral theory of block Jacobi matrices, then introduce the various distributions which we will work with.

Throughout the paper, we will use $\FF$ to denote $\R$ ($\beta=1$) or $\CC$ ($\beta=2$). In particular, we use $\FF$-Hermitian and $\FF$-unitary for real symmetric/Hermitian and orthogonal/unitary matrices. We use $\mbf{X}^\T$ to denote the transpose/conjugate transpose of an $\FF$-matrix $\mbf{X}$.

\subsection{Block Jacobi matrices}

We work with the following block generalization of tridiagonal Jacobi matrices.

\begin{Definition}\label{def:tridiag}
 Let $r, n\ge 1$. An $(rn)\times(rn)$ matrix $\mbf{T}$ is called an $\FF$-valued $r$-block Jacobi matrix if it is a $\FF$-Hermitian block tridiagonal matrix built from $r\times r$ blocks satisfying the following conditions. The diagonal blocks $\mbf{A}_1, \dots, \mbf{A}_n$ are $r\times r$ $\FF$-Hermitian matrices. The off-diagonal blocks $\mbf{B}_1, \dots, \mbf{B}_{n-1}$ above the diagonal are lower triangular with positive diagonal entries, see \eqref{eq:T}.
 We denote the set of such matrices by $\mathfrak{M}_{n,\beta, r}$,
 \begin{align}\label{eq:T}
 \mbf{T}=
 \left[\begin{array}{ccccc}
 \mbf{A}_1& \mbf{B}_1 & 0 &\dots & \\
 \mbf{B}_1^{\dag} & \mbf{A}_2 &\mbf{B}_2 &\dots \\
 0&\ddots & \ddots & \ddots &0 \\
 & 0 & \mbf{B}_{n-2}^\dag &\mbf{A}_{n-1} &\mbf{B}_{n-1} \\
 & & 0 & \mbf{B}_{n-1}^\dag & \mbf{A}_n\\
 \end{array} \right].
 \end{align}
 \end{Definition}
Note that an $r$-block Jacobi matrix can be viewed $(2r+1)$-diagonal band matrix with positive entries at the boundaries of the band.

Let \smash{$\mbf{e}_{\lst{r}}=[\mbf{I}_r,\mbf{0}_{r\times (n-1)r}]^{\T}$} denote $(rn)\times r$ matrix built from the first $r$ coordinate vectors. (We do not explicitly denote the $n$-dependence.)

The proof of the following theorem can be found, for example, in \cite{Spike2}, it relies on the Householder tridiagonalization algorithm in a block setting.
\begin{Theorem}[\cite{Spike2}]\label{thm:block_basic_1}
 Suppose that $\mbf{M}$ is an $\FF$-Hermitian $rn\times rn$ matrix for which the matrix
 \begin{align}\label{eq:S1234}
 \mbf{S}=\big[\mbf{e}_{\lst{r}}, \mbf{M}\mbf{e}_{\lst{r}},\dots, \mbf{M}^{n-1}\mbf{e}_{\lst{r}}\big]
 \end{align}
 is invertible.
 Then there is an $\FF$-unitary matrix $\mbf{O}$ of the form $\mbf{I}_r\oplus \widetilde{\mbf{O}}$ and a unique $\mbf{T}\in \mathfrak{M}_{n,\beta, r}$, so that $\mbf{T}=\mbf{O}^{\T} \mbf{M} \mbf{O}$. The matrix $\mbf{O}$ can be chosen as the $\mbf{Q}$ in
 the unique QR decomposition~${\mbf{S}=\mbf{Q}\mbf{R}}$ for which $\mbf{R}$ has positive diagonal entries.
\end{Theorem}

For $r=1$, the spectral measure of an $n\times n$ tridiagonal Hermitian matrix $\mbf{T}$ with respect to the first coordinate vector $\mbf{e}_1$ is defined as the probability measure
\begin{align}\label{eq:spec_m}
\mu=\sum_{j=1}^n |\mbf{v}_{j,1}|^2 \delta_{\lambda_j}.
\end{align}
Here $\mbf{v}_{j,1}$ is the first coordinate of the normalized eigenvector corresponding to $\lambda_j$. Our next definition provides a natural extension of the spectral measure for $r$-block Jacobi matrices.

\begin{Definition}
 Suppose that $\mbf{M}$ is an $\FF$-Hermitian $rn\times rn$ matrix.
 We define the spectral measure of $\mbf{M}$ with respect to $\mbf{e}_{\lst{r}}$ as the $r\times r$ matrix-valued measure
\begin{align*}
 \mu_{\lst{r}}=\sum_{j=1}^{rn} \mbf{v}_{j,\lst{r}} \cdot \mbf{v}_{j,\lst{r}}^{\T} \delta_{\lambda_j}.
\end{align*}
Here $\mbf{v}_{j}$ is the normalized eigenvector corresponding to $\lambda_j$, and
$\mbf{v}_{j,\lst{r}}\in \FF^r$ is the projection of~$\mbf{v}_j$ to the first $r$ coordinates.
\end{Definition}
Note that $\mu_{\lst{r}}$ only depends on the eigenspaces, so it is well-defined even though the choice of $\mbf{v}$ is not unique. If $\mbf{T}$ is the $r$-block Jacobi matrix obtained from an $\FF$-Hermitian $\mbf{M}$ via Theorem~\ref{thm:block_basic_1} then we have
\begin{align*}
\int x^j {\rm d}\mu_{\lst{r}}=\mbf{e}_{\lst{r}}^{\T} \mbf{M}^j \mbf{e}_{\lst{r}}= \mbf{e}_{\lst{r}}^{\T} \mbf{T}^j \mbf{e}_{\lst{r}}.
\end{align*}
It can be shown that there is a one-to-one correspondence between the $r$-block Jacobi matrices and possible $r\times r$ matrix valued `probability' measures, see \cite[Section 2]{MOPUC}.

\subsection{Random block matrices}\label{subs:matrix_distr}

We start with an overview of the various distributions that serve as building blocks for our models, and then provide a precise definition of the $\HH(r,s)$ and $\WW(r,s)$ distributions.

\begin{Definition}
 The $\FF$-valued standard normal is denoted by $\FF N(0,1)$. The components are independent mean zero normals with variance $\frac{1}{\beta}$. The probability density function is proportional to~\smash{${\rm e}^{-\frac{\beta}{2} |x|^2}$}.
\end{Definition}

We record the fact that if $\mbf{x}$ is a $d$-dimensional random vector with i.i.d.~$\FF N(0,1)$ entries then the distribution of $|\mbf{x}|$ is \smash{$\frac{1}{\sqrt{\beta}}\chi_{\beta d}$}. The probability density function of $|\mbf{x}|$ is
 \[
2 \frac{ (\beta/2)^{\frac{\beta d}{2}}}{\Gamma(\beta d/2)} x^{\beta d-1} {\rm e}^{-\frac{\beta}{2} x^2}.
 \]

\begin{Definition}
 Let $\mbf{Y}$ be an $n\times n$ matrix with i.i.d.~$\FF N(0,1)$ entries, and set \smash{$\mbf{X}=\frac1{\sqrt{2}} \bigl(\mbf{Y}+\mbf{Y}^{\T}\bigr)$}. The distribution of $\mbf{X}$ is called the $\FF$-valued Gaussian ensemble, or G$\FF$E$(n)$. For $\beta=1$, this is the Gaussian orthogonal ensemble (GOE), and for $\beta=2$, this is the Gaussian unitary ensemble (GOE).
\end{Definition}

The diagonal entries of G$\FF$E are $N\bigl(0,\frac{2}{\beta}\bigr)$ distributed, while the off-diagonal entries are i.i.d.~$\FF N(0,1)$. The entries are independent up to the real/Hermitian symmetry. In the matrix variables, the probability density function of G$\FF$E is proportional to \smash{$ {\rm e}^{-\frac{\beta}{4} \Tr \mbf{X}\mbf{X}^{\T}}$}.

\begin{Definition}
 Let $\mbf{Y}$ be an $n\times m$ (with $n\le m$) matrix with i.i.d.~$\FF N(0,1)$ entries. The distribution of the matrix $\mbf{X}=\mbf{Y}\mbf{Y}^{\mathsf T}$ is called the $\FF$-valued Wishart distribution with parameters~$(n,m)$.
\end{Definition}

The following is a classical result in random matrix theory.
\begin{Theorem}
 The joint eigenvalue density of the $\FF$-valued $n\times n$ Gaussian ensemble is given by \eqref{eig_DE}. The distribution is called the Gaussian beta ensemble, and it is denoted by $G{\beta}E(n)$.

The joint eigenvalue density of the $\FF$-valued Wishart distribution with parameters $(n,m)$ is given by \eqref{eig_DE1}.
 The distribution is called the Laguerre beta ensemble, and it is denoted by~$L{\beta}E(n,m)$.

In both cases, the normalized eigenvectors can be chosen in a way so that the eigenvector matrix is Haar-distributed on the $n\times n$ $\FF$-unitary matrices while being independent of the eigenvalues.
\end{Theorem}

\begin{Definition}
The $\FF$-valued square root Wishart matrix with parameters $ m> n-1$ is the distribution of the $n\times n$ lower triangular matrix $\mbf{X}$ with the following independent entries:
\begin{align*}
 x_{i,j}\sim \begin{cases}
 \FF N(0,1)& \text{if $i>j$},\\
 \dfrac{1}{\sqrt{\beta}} \chi_{\beta (m+1-i)}& \text{if $i=j$},\\
 0& \text{if $i<j$}.
 \end{cases}
\end{align*}
 We denote this distribution by $\SQW(n,m)$.
\end{Definition}

We note that the joint probability density function of the non-zero entries of $\SQW(n,m)$ is proportional to
\begin{align*}
 \prod_{i>j} {\rm e}^{-\frac{\beta}{2} |x_{i,j}|^2} \prod_{i=1}^n x_{i,i}^{\beta (m+1-i)-1} {\rm e}^{-\frac{\beta}{2} x_{i,i}^2}={\rm e}^{-\frac{\beta}{2} \Tr \mbf{X}\mbf{X}^\T} \det(\mbf{X})^{\beta (m+1)-1} \prod_{i=1}^n x_{i,i}^{-\beta i}.
\end{align*}

As the following classical result due to Bartlett \cite{Bartlett1933} shows, $\SQW(n,m)$ gives the distribution of the Cholesky factor of the Wishart distribution.
\begin{Theorem}[\cite{Bartlett1933}]\label{thm:bartlett}
 Suppose that the matrix $\mbf{X}$ has $\FF$-valued Wishart distribution with parameters $(n,m)$. Let $\mbf{R}$ be the lower triangular square root of $\mbf{X}$ with almost surely positive diagonal entries: $\mbf{X}=\mbf{R} \mbf{R}^{\T}$. Then $\mbf{R}$ has $\SQW(n,m)$ distribution.
\end{Theorem}

We can now introduce the family of random block matrices that we study.

\begin{Definition}
\label{def:BlockH}
Let $r,n\ge 1$ and $s\ge 0$.
We denote by $\HH(r,s)$ the distribution of the $\FF$-valued random $r$-block Jacobi matrix of size
$(rn)\times(rn)$ with independent blocks $\mbf{A}_k$,
 $\mbf{B}_k$ where~${\mbf{A}_k\sim}$ G$\FF$E$(r)$ and $\mbf{B}_k\sim \SQW(r,(r+s)(n-k))$.
\end{Definition}

Note that $\HH(1,0)$ is just the distribution of the tridiagonal matrix of Dumitriu and Edelman (and Trotter) given for the Gaussian beta ensemble. As the following theorem shows, for $r\ge 1$ the $\HH(r,0)$ distribution is the result of the $r$-block Householder process applied to G$\FF$E$(rn)$.

\begin{Theorem}[\cite{Spike2}]\label{thm:GFE_block}
Let $\mbf{M}$ have G$\FF$E$(rn)$ distribution, and consider the matrix $\mbf{S}$ defined via~\eqref{eq:S1234}. Then $\mbf{S}$ is a.s.~invertible, and the $r$-block Jacobi matrix $\mbf{T}$ produced by Theorem~{\rm\ref{thm:block_basic_1}} has~$\HH(r,0)$ distribution.

The eigenvalues of $\mbf{T}$ are distributed as $G\beta E(rn)$, and the normalized eigenvector matrix $\mbf{V}=[\mbf{v}_{i,j}]_{i,j\in \lst{rn}}$ can be chosen in a way so that the first $r$ rows of $\mbf{V}$ are independent of the eigenvalues and have the same distribution as the first $r$ rows of an $rn\times rn$ Haar $\FF$-unitary matrix.
\end{Theorem}

Theorem~\ref{thm:GFE_block} fully describes the distribution of the matrix valued spectral measure $\mu_{\lst{r}}$ of~$\mbf{T}$. In particular, it shows that the weights and the support are independent of each other, and the weights can be obtained from a Haar $\FF$-unitary matrix.

\begin{Definition}\label{def:BlockW}
Let $r,n\ge 1$, $m>n-1/r$, and $s\ge 0$.
Let $\mathbf{L}$
be an $rn \times rn$ block bidiagonal matrix with independent $r\times r$ diagonal and upper offdiagonal blocks denoted by
$\{\mbf{D}_i\}_{i=1,n}$ and~${\{\mbf{O}_i\}_{i=1, n-1}}$ with \smash{$\mbf{D}_i^{\T}\sim \SQW(r,(r+s)(m+1-i))$} and $\mbf{O}_i\sim \SQW(r,(r+s)(n-i))$. We denote the distribution of $\mbf{W}=\mbf{L}\mbf{L}^{\T}$ by $\WW(r,s)$.
\end{Definition}
Again, $\WW(1,0)$ is just the tridiagonal model given by Dumitriu and Edelman for the Laguerre beta ensemble.
The analogue of Theorem~\ref{thm:GFE_block} holds.
\begin{Theorem}[\cite{Spike2}]\label{thm:W_block}
 Let $\mbf{M}$ have $\FF$-valued Wishart distribution with parameters $(rn,rm)$, and consider the matrix $\mbf{S}$ defined via \eqref{eq:S1234}. Then $\mbf{S}$ is a.s.~invertible, and the $r$-block Jacobi matrix~$\mbf{T}$ produced by Theorem~{\rm\ref{thm:block_basic_1}} has $\WW(r,0)$ distribution.
 The eigenvalues of $\mbf{T}$ are distributed as~$L\beta E(rn,rm)$, and the normalized eigenvectors can be chosen in a way that the first $r$ rows are independent of the eigenvalues and have the same distribution as the first $r$ rows of an $rn\times rn$ Haar $\FF$-unitary matrix.
\end{Theorem}

\section{New distributions via biasing}\label{sec3}

We start this section with a brief review of the Dumitriu--Edelman result \cite{DE}. We introduce the key tools for our block generalization and provide the proofs of our main theorems modulo a~certain moment computation that is delayed to the subsequent sections.

\subsection{Revisiting the Hermite beta ensemble}

For completeness, we state the Dumitriu--Edelman result in full and provide a proof which foreshadows the techniques used to prove Theorem~\ref{thm:main}.

\begin{Theorem}[\cite{DE}]\label{thm:DE}
Fix $\beta>0$ and an integer $n\ge 1$. Let $a_1,\dots, a_n$, $b_1, \dots, b_{n-1}$ be independent random variables with $a_j\sim N\bigl(0,\frac{2}{\beta}\bigr)$, \smash{$b_j\sim \frac{1}{\sqrt{\beta}}\chi_{\beta (n-j)}$}. Then the symmetric tridiagonal matrix $\mbf{T}$ with diagonal $a_1,a_2,\dots$ and off-diagonal $b_1,b_2, \dots$ has a joint symmetrized eigenvalue density on $\R^n$ given by
\begin{align*}%\label{eq:GbE}
 \frac{1}{Z_{n,\beta}} \left|\Delta(\lambda)\right|^\beta {\rm e}^{-\frac{\beta}{4} \sum_{j=1}^n \lambda_j^2},
\end{align*}
with
\begin{align}\label{eq:GbE_constant}
 Z_{n,\beta}={n!} (2\pi)^{n/2} (\beta/2)^{-\frac{\beta}{4}n(n-1)-\frac{n}{2}} \Gamma(\beta/2)^{-n} \prod_{j=1}^n \Gamma(\beta j/2).
\end{align}
Moreover, the spectral weights of $\mbf{T}$ corresponding to the first coordinate vector have joint distribution given by $\mathrm{Dirichlet}(\beta/2,\dots, \beta/2)$, and this weight vector is independent of the eigenvalues.
\end{Theorem}

\begin{proof}
Consider an $n\times n$ Jacobi matrix $\mbf{T}$
with diagonal entries $a_1,\dots, a_n$ and off-diagonal positive entries $b_1, \dots, b_{n-1}$. Denote by $p_j$ the spectral weight of $\lambda_j$ in the spectral measure~\eqref{eq:spec_m}. It is well known that
\begin{align}\label{eq:magic_Delta_p}
|\Delta({\lambda})|= \prod_{k=1}^n p_k^{-1/2} \prod_{k=1}^{n-1} b_k^{(n-k)},
\end{align}
see, for instance, \cite[equation~(1.148)]{ForBook}.
We also take as given that the theorem holds for $\beta=1$ due to \cite{Trotter}, and the fact that the
Householder tridiagonalization process does not change the spectral measure with respect to the first coordinate.

Next, for $\mbf{T}$ a random tridiagonal matrix defined in the statement with $\beta=1$, introduce a~biased version of the distribution of $\mbf{T}$ with the biasing function
\[
g_\beta(\mbf{b})=\prod_{k=1}^{n-1} b_k^{(\beta-1)(n-k)}.
\]
The biasing produces a random tridiagonal matrix $\mbf{\wt{T}}$ where the diagonal and off-diagonal entries are still independent, the distribution of the diagonal entries is still $N(0,2)$, but the distribution of the $k$th off-diagonal entry has changed from $\chi_{n-k}$ to $\chi_{\beta(n-k)}$. By \eqref{eq:magic_Delta_p}, we have
\begin{align}\label{eq:bias_DE}
g_\beta(\mbf{b})=|\Delta({\lambda})|^{\beta-1} \prod_{k=1}^n p_k^{-\frac{\beta-1}{2}},
\end{align}
hence biasing the entries of $\mbf{T}$ with $g_\beta(\mbf{b})$ is the same as biasing the spectral variables $\lambda, \mbf{p}$ with the appropriate product on the right-hand side of \eqref{eq:bias_DE}. This immediately implies that the vector of eigenvalues and the vector of spectral weights of $\mbf{\wt{T}}$ are still independent of each other, that the joint eigenvalue density of $\mbf{\wt{T}}$ is proportional to $|\Delta(\lambda)|^\beta {\rm e}^{-\frac{1}{4}\sum_{k=1}^n \lambda_k^2}$, and that its spectral weights have Dirichlet$(\beta/2,\dots,\beta/2)$ distribution.

The complete statement of the theorem now follows after scaling $\mbf{\wt{T}}$ by \smash{$ \frac{1}{\sqrt{\beta}}$}. The value of the normalizing constant $Z_{n,\beta}$ follows from the known $\beta=1$ factor (see \cite[equation~(1.160)]{ForBook}) along with an evaluation of $E[g_\beta(\mbf{b})]$.
\end{proof}

\begin{Remark}
The original proof of Theorem~\ref{thm:DE} given in \cite{DE} is slightly different. It uses the~${\beta=1}$ case as a starting point to derive an expression for the Jacobian of the one-to-one transformation~${(\mbf{a},\mbf{b}) \mapsto ({\lambda}, \mbf{p})}$. It then uses the identity \eqref{eq:magic_Delta_p} to
compute the joint density of the spectral variables $({\lambda}, \mbf{p})$ for the tridiagonal matrix $\mbf{T}$. As it was already remarked in \cite{DE}, the expression for the Jacobian of the transformation $(\mbf{a},\mbf{b}) \mapsto ({\lambda}, \mbf{p})$
can also be computed directly, without relying on Theorem~\ref{thm:DE} being true in the $\beta=1$ case.
\end{Remark}

\subsection{Key spectral identity}

We establish a block Jacobi matrix equivalent of the identity \eqref{eq:magic_Delta_p} which relates products of powers of the determinants of the off-diagonal blocks to eigenvalues and eigenvector entries. First, we need the following definition.

\begin{Definition}
For $\lambda = \{ \lambda_j\}_{j \in \lst{rn} }$ and an $r\times rn $ matrix $\mbf{X}$, define the $rn \times rn$ matrix $\mbf{M}=\mbf{M}({\lambda}, \mbf{X})$ entry-wise as
\begin{align}\label{eq:MagicM_ij}
 M_{i,j}=\lambda_i^{\nint{j}{r}} x_{ \nfr{j}{r},i},
\end{align}
in which
\begin{align*}
 % \lfloor i \rfloor_k
 \nint{j}{r}
 :=\left\lfloor \frac{j-1}{r} \right\rfloor, \qquad \nfr{j}{r}:=j-r \nint{j}{r}.
\end{align*}
Note that the shift in $ \nint{j}{r}$ means that $\lfloor r \rfloor_r = 0$ and hence
$\nfr{j}{r}\in \lst{r}$.
For an $rn \times rn$ matrix~$\mbf{X}$, we define $\mbf{M}({\lambda}, \mbf{X})$ using the $r \times rn$ submatrix of the first $r$ rows of $\mbf{X}$.
\end{Definition}

The operation $({\lambda}, \mbf{X})\mapsto \mbf{M}$ can be realized as
 a block-matrix built from the row-by-row Kronecker products of the $rn\times r$ matrix $\mbf{X}^{\mathsf T}$ and the
 $rn\times n$ Vandermonde type matrix $\mbf{\Lambda}=\smash{\bigl(\lambda_{i}^{j-1}\bigr){}_{i\in \lst{rn},j\in \lst{n}}}$. Similar constructions show up in the literature as the ``face-splitting product'' or Khatri--Rao product~\cite{KR_1968}.

As an example,
\begin{align*}%\label{eq:M6}
 \mbf{M}({\lambda}, \mbf{X})=
 \left(
\begin{array}{cccccc}
 x_{1,1} & x_{2,1} & \lambda_1 x_{1,1} & \lambda_1 x_{2,1} & \lambda_1^2 x_{1,1} & \lambda_1^2 x_{2,1} \\
 x_{1,2} & x_{2,2} & \lambda_2 x_{1,2} & \lambda_2 x_{2,2} & \lambda_2^2 x_{1,2} & \lambda_2^2 x_{2,2} \\
 x_{1,3} & x_{2,3} & \lambda_3 x_{1,3} & \lambda_3 x_{2,3} & \lambda_3^2 x_{1,3} & \lambda_3^2 x_{2,3} \\
 x_{1,4} & x_{2,4} & \lambda_4 x_{1,4} & \lambda_4 x_{2,4} & \lambda_4^2 x_{1,4} & \lambda_4^2 x_{2,4} \\
 x_{1,5} & x_{2,5} & \lambda_5 x_{1,5} & \lambda_5 x_{2,5} & \lambda_5^2 x_{1,5} & \lambda_5^2 x_{2,5} \\
 x_{1,6} & x_{2,6} & \lambda_6 x_{1,6} & \lambda_6 x_{2,6} & \lambda_6^2 x_{1,6} & \lambda_6^2 x_{2,6} \\
\end{array}
\right),
\end{align*}
illustrates the $n=3$, $r=2$ case. The advertised formula, recorded next, involves the determinant of $\mbf{M}$.

\begin{Proposition}\label{prop:magic}
Suppose that $\mbf{T}\in \mathfrak{M}_{n,\beta, r}$ with blocks $\mbf{A}_j$, $\mbf{B}_j$, recall Definition {\rm\ref{def:tridiag}}. Then
\begin{align}\label{eq:magic11}
 \prod_{j=1}^{n-1} \det(\mbf{B}_j)^{n-j}=|{\det \mbf{M}({\lambda}, \mbf{Q})}|,
\end{align}
where $\mbf{\lambda}$ are the eigenvalues of $\mbf{T}$, and $\mbf{Q}$ is the matrix of the normalized eigenvectors of $\mbf{T}$ ordered according to $\lambda$.
\end{Proposition}

When $r=1$, one can easily see that $\det {\mbf{M}(\lambda, \mbf{Q})} = \Delta(\lambda) \prod_{i=1}^n q_{1,i}$ upon factoring out the eigenvector coordinates from each row, recovering the identity \eqref{eq:magic_Delta_p}.

We remark that for a particular $\mbf{T}$ the choice of $\mbf{Q}$ is not unique: eigenvectors could be multiplied by a unit phase, and there is even more freedom if $\mbf{T}$ has eigenvalues with multiplicity higher than one. The proof below shows that the expression on the right-hand side of \eqref{eq:magic11} does not depend on the particular choice of $\mbf{Q}$.

\begin{proof}[Proof of Proposition~\ref{prop:magic}]
 Consider the $rn \times rn $ matrix $\mbf{S}$ from \eqref{eq:S1234} built from $\mbf{T}$. This is an upper triangular block matrix where the $j$th diagonal block is $\mbf{B}_{j-1}^{\T}\cdots \mbf{B}_1^{\T}$, the first diagonal block being $\mbf{I}_r$. The determinant of $\mbf{S}$ is the product of the determinants of its diagonal blocks
 \[
 \det \mbf{S}=\prod_{j=2}^n \det\bigl(\mbf{B}_{j-1}^{\T}\cdots \mbf{B}_1^{\T}\bigr)=\prod_{j=1}^{n-1} \det(\mbf{B}_j)^{n-j}.
 \]
Denote by $\boldsymbol{\Lambda}$ the diagonal matrix built from $\lambda$. Then we have $\mbf{T}=\mbf{Q}\boldsymbol{\Lambda} \mbf{Q}^{\T}$ and
 \[
 \mbf{T}^{j} \mbf{e}_{\lst{r}}=
 \mbf{Q} \boldsymbol{\Lambda}^{j} \mbf{Q}^{\T}\mbf{e}_{\lst{r}}.
 \]
Hence $\mbf{S}= {\mbf{Q}} \widetilde{\mbf{S}}$, and \smash{$|{\det \mbf{S}}|=\bigl|{\det \wt{\mbf{S}}}\bigr|$} with
\begin{equation*}
 \widetilde{\mbf{S}}=\big[\mbf{Q}^{\T}\mbf{e}_{\lst{r}}, \boldsymbol{\Lambda}\mbf{Q}^{\T}\mbf{e}_{\lst{r}}, \dots, \boldsymbol{\Lambda}^{n-1}\mbf{Q}^{\T}\mbf{e}_{\lst{r}}\big]=\mbf{M}\bigl(\mbf{\lambda}, \overline{\mbf{Q}}\bigr).
\end{equation*}
Since $\big|{\det \mbf{M}\big(\mbf{\lambda}, \overline{\mbf{Q}}\big)}\big|=\big|{\det \mbf{M}(\mbf{\lambda}, {\overline{}\mbf{Q}})} |$, the proposition follows.
\end{proof}

We remark that for $r=2$ and $\lambda_j>0$ the following remarkable identity holds for $\det \mbf{M}({\lambda}, \mbf{Q})$ (see \cite[equation (52)]{LT2006})
\begin{align}\label{eq:Pfaffian_detM}
 \det \mbf{M}({\lambda}, \mbf{Q})=\prod_{1\le i<j\le 2n}\bigl(\sqrt{\lambda_i}+\sqrt{\lambda_j}\bigr) \operatorname{Pf}\left(\frac{q_{1,i}q_{2,j}-q_{1,j}q_{2,i}}{\sqrt{\lambda_i}+\sqrt{\lambda_j}}\right).
\end{align}
We are not aware of extensions of this result for $r>2$.

\subsection[Proofs of Theorem 1.1 and 1.6]{Proofs of Theorem~\ref{thm:main} and \ref{thm:main_W}}

We are now ready to present the proofs of our main theorems, modulo one key ingredient.

\begin{proof}[Proof of Theorem~\ref{thm:main} -- first steps]
The joint density of the entries of a matrix $\mbf{T}$ distributed as $\HH(r,s)$ is given by
\begin{align*}
 f(\mbf{A},\mbf{B})=C_{\beta,n,r,s} \exp\left(-\frac{\beta}{4} \Tr \mbf{T} \mbf{T}^{\T}\right) \prod_{m=1}^{n-1} \bigl(\det \mbf{B}_m\bigr)^{\beta(r+s)(n-m)+\beta-1} \prod_{m=1}^{n-1} \prod_{j=1}^r \bigl(b_{j,j}^{(m)}\bigr)^{-\beta j},
\end{align*}
where $C_{\beta, n, r, s}$ is an explicitly computable normalizing constant.

Note that the distribution $\HH(r,s)$ can be realized as the biased version of the distribution~$\HH(r,0)$ with the biasing function
\smash{$
\prod_{m=1}^{n-1} \bigl(\det \mbf{B}_m\bigr)^{\beta s(n-m)}$}.
Using $E_{s}$ for the expectation with respect to $\HH(r,s)$ (suppressing the dependence on $\beta$, $n$, $r$), we have, for any (bounded) test function $h$,
\begin{align*}
E_{s}[h(\mbf{T})] & =\frac{E_{0}\big[h(\mbf{T})\prod_{m=1}^{n-1} \bigl(\det \mbf{B}_m\bigr)^{\beta s(n-m)} \big]} {E_{0}\big[\prod_{m=1}^{n-1} \bigl(\det \mbf{B}_m\bigr)^{\beta s(n-m)} \big]}
 := \frac{1}{Z_n} E_{0}\left[h(\mbf{T})\prod_{m=1}^{n-1} \bigl(\det \mbf{B}_m\bigr)^{\beta s(n-m)} \right].
\end{align*}
In particular, this is true if $h(\mbf{T})=g({\lambda})$ with a (bounded) test function $g$. By Proposition~\ref{prop:magic}, we have
\begin{align}\label{eq:magic_identity}
 \prod_{m=1}^{n-1} \bigl(\det \mbf{B}_m\bigr)^{\beta s(n-m)}=|{\det \mbf{M}({\lambda}, \mbf{Q})}|^{\beta s},
\end{align}
where $\mbf{Q}$ is the $r\times (rn)$ matrix of the first $r$ coordinates of the normalized eigenvectors of $\mbf{T}$. By Theorem~\ref{thm:GFE_block}, $\mbf{Q}$ and ${\lambda}$ are independent under the distribution $\HH(r,0)$, and $\mbf{Q}$ is distributed as the first $r$ rows of a Haar $\FF$-unitary matrix. This implies that
\begin{align*}
 E_{s}[g({\lambda})] & = \frac{1}{Z_n} E_{0}\left[g({\lambda}) F_{\beta,r,n, \beta s}({\lambda})\right],
\end{align*}
where
 \begin{align}
 \label{Reweight}
 F_{\beta, r,n, \beta s}({\lambda}) & =
 E_{\boldsymbol{Q}} \big[|{\det \mbf{M}({\lambda}, \mbf{Q})}|^{\beta s}\big].
\end{align}
In other words, the distribution of ${\lambda}$ is just G$\beta$E biased by the functional $
F_{\beta, r,n,\beta s}({\lambda})$.

As for the normalizer, the left-hand side of \eqref{eq:magic_identity} readily yields
\begin{align}
 Z_n &= E_{0} \left[ \prod_{m=1}^{n-1}\bigl(\det \mbf{B}_m\bigr)^{(s)(n-m)} \right]
 = \prod_{m=1}^{n-1} \prod_{i=1}^r E \bigl(b^{(m)}_{i,i}\bigr)^{\beta s(n-m)}\nonumber \\
 & =
\prod_{m=1}^{n-1} \prod_{i=1}^r
(\beta/2) ^{-\frac{\beta s}{2}(n-m)}\frac{ \Gamma \bigl(\frac{\beta}{2}((r+s)(n-m)-i+1 ))\bigr)}{\Gamma \bigl(\frac{\beta}{2}{((r(n-m)-i+1) }\bigr)} \nonumber
\\&=(\beta/2)^{-\frac{\beta }{4} r s n(n-1)} \prod_{m=1}^{n-1} \prod_{i=1}^r
\frac{ \Gamma \bigl(\frac{\beta}{2}((r+s)(n-m)-i+1 )\bigr)}{\Gamma \bigl(\frac{\beta}{2}(r(n-m)-i+1) \bigr)}
\label{eq:bias_const}
\end{align}
since the $b^{(m)}_{i,i}$ are independent \smash{$ \frac{1}{\sqrt{\beta}} \times \chi_{\beta (r(n-m)+1-i)}$}
variables.
\end{proof}

Identifying the joint eigenvalue density for any $\HH(r,s)$ model then comes down to computing the determinant moments indicated in \eqref{Reweight}. We have been able to do this only for general~${r\ge 2}$ and for $\beta s=2$, and $r=2$, $\beta s=4$.
In particular, the proof of Theorem~\ref{thm:main} is finished by establishing the following.

\begin{Theorem}\label{thm:moment}
With $\boldsymbol{Q}$ Haar distributed $\FF$-unitary matrix, we have
\begin{align*}
 E_{\mbf{Q}} | \det \mbf{M}({\lambda}, \mbf{Q}) |^{\beta s} = c_{n, \beta, r,\beta s} \times
 \begin{cases}
 \sum\limits_{(\cA_1,\dots, \cA_r)\in \cP_{r,n}} \prod_{j=1}^r \Delta(\cA_j)^2 & \mbox{for } r \ge2, \ \beta s=2,\\[10pt]
 \sum\limits_{(\cA,\cA')\in \cP_{2,n}} \Delta(\cA)^4 \Delta(\cA')^4 & \mbox{for } r=2, \ \beta s = 4,
 \end{cases}
\end{align*}
where
\begin{align*}
 c_{n,\beta, r, \beta s} =(\beta/2) ^{\frac{\beta}{2}rsn} \prod_{i=1}^{r}
\frac{\Gamma \bigl(\frac{\beta }{2}(rn+1-i)\bigr)}{ \Gamma \bigl(\frac{\beta}{2}((r+s)n+1-i) \bigr)} \times \begin{cases}
 1, &r\ge2,\ \beta s=2,\\
 (12/\beta)^n, &r=2,\ \beta s=4.
\end{cases}
\end{align*}

\end{Theorem}
One immediately recognizes the interaction term reported in \eqref{density1}. To arrive at the form of~${r=2}$, $s=2,4$ densities appearing in \eqref{density2}, we will also prove
\begin{equation}
\label{id_quad_to_square}
 \sum_{ (\cA,\cA')\in \cP_{2,n}
 } \Delta (\cA)^4 \Delta (\cA')^4
 = 2^{-n} \left(\sum_{ (\cA,\cA')\in \cP_{2,n}
 }
 \Delta (\cA)^2 \Delta (\cA')^2\right)^2,
\end{equation}
and
 \begin{align}
 \label{id_pfaff}
 \sum_{
 (\cA,\cA')\in \cP_{2,n}
 } \Delta(\cA)^2 \Delta (\cA')^2&=(-2)^n \Delta (\lambda) \Pf\left(\frac{\ind_{i\neq j}}{\lambda_i-\lambda_j}\right).
 \end{align}
The proof of Theorem~\ref{thm:moment} is divided between Sections~\ref{sec:2moment} and~\ref{sec:4moment}, for $s=2$ and $s=4$, respectively. The identities \eqref{id_quad_to_square} and \eqref{id_pfaff} are proven in Section \ref{sec:det_identities}, see in particular Lemma~\ref{lem:det4_identities} along with the related
Propositions~\ref{prop:det2+Pf} and~\ref{prop:det4+det2}.

\begin{proof}[Proof of Theorem~\ref{thm:main} -- final steps]
We have shown that the joint eigenvalue distribution of $\HH(r,s)$ is just $G\beta{E}(rn)$ biased by the functional $
F_{\beta, r,n,\beta s}({\lambda})$ given in \eqref{Reweight}. Theorem~\ref{thm:moment} provides $
F_{\beta, r,n,\beta s}({\lambda})$ for $r\ge 2$ and $\beta s=2$ and for $r=2$, $\beta s=4$.

When $r\ge 2$, $\beta s=2$, this gives the joint density function \eqref{density1}. The normalizing constant~$Z_{n,\beta, r, \beta s}$ satisfies
\begin{align}\label{eq:norm_const}
 Z_{n,\beta,r,\beta s}=Z_{nr, \beta} Z_n c_{n,\beta, r, \beta s}^{-1},
\end{align}
where $Z_{nr,\beta}$ is the normalizing constant in \eqref{eq:GbE_constant}, $Z_n$ is given in \eqref{eq:bias_const} and $c_{n,\beta, r, \beta s}$ is given in Theorem~\ref{thm:moment}.

When $r=2$ and $\beta s=2$ or 4, then we also use \eqref{id_quad_to_square} and \eqref{id_pfaff} to rewrite expression in Theorem~\ref{thm:moment}
with a Pfaffian to obtain the joint density function \eqref{density2}. The normalizing constant is again given by \eqref{eq:norm_const}, note that we needed the additional $2^n$ factor in \eqref{density2} to match the two forms of the $r=2$, $\beta s=2$ case. The reported constant in Theorem~\ref{thm:main} now follows after some algebra, noting that
\begin{gather*}
 \prod_{m=1}^{n-1}\prod_{i=1}^r
\frac{ \Gamma \bigl(\frac{\beta}{2}((r+s)(n-m)-i+1 )\bigr)}{\Gamma \bigl(\frac{\beta}{2}(r(n-m)-i+1) \bigr)} \prod_{i=1}^{r}
\frac{ \Gamma \bigl(\frac{\beta}{2}((r+s)n+1-i) \bigr)}{\Gamma \bigl(\frac{\beta }{2}(rn+1-i)\bigr)}
\\
\qquad = \prod_{k=1}^{rn} \frac{\Gamma\bigl(\frac{\beta}{2} \bigl(k+s \lceil\frac{k}{r} \rceil\bigr)\bigr)}{\Gamma\bigl(\frac{\beta}{2} k\bigr)}.\tag*{\qed}
\end{gather*}\renewcommand{\qed}{}
 \end{proof}

The proof of Theorem~\ref{thm:main_W} follows the same biasing idea as the proof of Theorem~\ref{thm:main}.

\begin{proof}[Proof of Theorem~\ref{thm:main_W}]
Let $\mbf{T}=\mbf{L}\mbf{L}^{\T}\sim \WW(r,s)$, with $\mbf{L}$ as in Definition \ref{def:BlockW}. (Note that~$\mbf{L}$ is a function of $\mbf{T}$.) The joint density of the entries $\mbf{O}_m$, $\mbf{D}_m$ of $\mbf{L}$ are given by
\begin{align*}
 f(\mbf{D}, \mbf{O})={}&C_{\beta, n,m,r,s} \exp\left(-\frac{\beta}{2}\Tr \mbf{T} \mbf{T}^{\T}\right)\prod_{j=1}^n \bigl(\det \mbf{D}_j\bigr)^{\beta (r+s)(m-j)+\beta-1} \\
 & \times \prod_{j=1}^{n-1} (\det \mbf{O}_j)^{\beta (r+s)(n-j)+\beta-1}\prod_{j=1}^n \prod_{k=1}^r \bigl(d^{(j)}_{k, k}\bigr)^{-\beta k} \prod_{j=1}^{n-1} \prod_{k=1}^r \bigl(o^{(j)}_{k, k}\bigr)^{-\beta k},
\end{align*}
where $C_{\beta,n,m,r,s}$ is an explicitly computable normalizing constant.

As before, we can realize $\WW(r,s)$ using a biased version of $\WW(r,0)$. Using $E_s$ for the expectation with respect to $\WW(r,s)$, we have, for any bounded test function $h$,
\begin{align*}
E_{s}[h(\mbf{T})] ={}&\frac{E_{0}\big[h(\mbf{T})\prod_{j=1}^n \bigl(\det \mbf{D}_j\bigr)^{\beta s(m-j)} \prod_{j=1}^{n-1} (\det \mbf{O}_j)^{\beta s(n-j)} \big]} {E_{0}\big[\prod_{j=1}^n \bigl(\det \mbf{D}_j\bigr)^{\beta s(m-j)} \prod_{j=1}^{n-1} (\det \mbf{O}_j)^{\beta s(n-j)} \big]} \\
:={}& \frac{1}{Z_{n,m}} E_{0}\left[h(\mbf{T})\prod_{j=1}^n \bigl(\det \mbf{D}_j\bigr)^{\beta s(m-j)} \prod_{j=1}^{n-1} (\det \mbf{O}_j)^{\beta s(n-j)} \right].
\end{align*}
The matrix $\mbf{T}=\mbf{L} \mbf{L}^{\T}$ has
off-diagonal blocks $\mbf{B}_j=\mbf{D}_j \mbf{O}_j$, hence by Proposition~\ref{prop:magic} we have
\begin{align*}
 \prod_{j=1}^{n-1} \det(\mbf{D}_j \mbf{O}_j)^{n-j}=\prod_{j=1}^{n} \det(\mbf{D}_j)^{n-j} \prod_{j=1}^{n-1} (\det \mbf{O}_j)^{n-j}=|{\det \mbf{M}({\lambda}, \mbf{Q})}|,
\end{align*}
where $\mbf{Q}$ is the matrix of the normalized eigenvalue matrix of $\mbf{T}$.
We also have
\begin{align*}
 \prod_{j=1}^{rn} \lambda_j=\det \mbf{T}=(\det \mbf{L})^2=\prod_{j=1}^n (\det \mbf{D}_j)^2.
\end{align*}
Using $h(\mbf{T})=g({\lambda})$ with a (bounded) test function $g$, we now get
\begin{align*}
 E_{s}[g({\lambda})] & = \frac{1}{Z_n} E_{0}\left[g({\lambda}) F_{\beta,r,n,m, s}({\lambda})\right],
\end{align*}
where
 \begin{align*}
 \label{Reweight1}
 F_{\beta, r,n, m,s}({\lambda}) & =\prod_{j=1}^{rn} \lambda_j^{\frac{\beta}{2}s(m-n)} \times
 E_{\boldsymbol{Q}} \big[|{\det \mbf{M}({\lambda}, \mbf{Q})}|^{\beta s} \big].
\end{align*}
Hence the distribution of $\lambda$ under $\WW(r,s)$ is just L$\beta$E$(rn,rm)$ biased by $F_{\beta, r,n, m,s}({\lambda})$. The statement of the theorem now follows from Theorem~\ref{thm:moment} and the identities \eqref{id_quad_to_square} and \eqref{id_pfaff}. The normalizing constant can be explicitly evaluated using a similar argument as in Theorem~\ref{thm:main}.
\end{proof}

\subsection{A reduction and a representation}

We end this section with two lemmas that play significant roles in our computations of the moments of $\det \mbf{M}({\lambda}, \mbf{Q})$. The first observation is that in computing moments of $\det \mbf{M}({\lambda}, \mbf{Q})$ we may replace the Haar distributed $\mbf{Q}$ with a matrix of independent Gaussians.

\begin{Lemma}[Gaussian reduction lemma]\label{lem:Gauss_red}
Let $\mbf{X}$ be an $rn \times rn$ matrix of $i.i.d.$ $\mathbb{F} N(0,1)$ entries, and $\mbf{Q}$ a Haar distributed $\mathbb{F}$-unitary matrix.
Let \smash{$R_{i,i}\sim \frac{1}{\sqrt{\beta}} \chi_{\beta(rn+1-i)}$}, $1\le i\le r$, be
independent of each other and of $\mbf{Q}$.
Then for any fixed collection $\lambda_i\in \R$, $i \in \lst{rn}$, we have the following distributional identity:
\begin{align*}
\det \mbf{M}({\lambda}, \mbf{X})\eqd \prod_{i=1}^{r} R_{i,i}^n \times
\det \mbf{M}({\lambda}, \mbf{Q}).
\end{align*}
It follows that for any $s>0$, we have
\begin{align*}
 E |{\det \mbf{M}({\lambda}, \mbf{Q})}|^{\beta s}=\kappa_{n,\beta, n,r, \beta s}
\cdot E |{\det \mbf{M}({\lambda}, \mbf{X})}|^{\beta s},
\end{align*}
with \smash{$\kappa_{n,\beta, n,r, \beta s} =(\beta/2) ^{\frac{1}{2}rsn}\times \prod_{i=1}^{r}
 \frac{\Gamma \left(\frac{\beta }{2}(rn+1-i)\right)}{ \Gamma \left(\frac{\beta}{2}((r+s)n+1-i) \right)}$}.
\end{Lemma}

\begin{proof}
Denote the columns of $\mbf{X}$ by $\mbf{X}_i$, $i\in \lst{rn}$. Consider the QR decomposition $\mathbf{X}=\mathbf{Q}\mathbf{R}$ of the matrix $\mathbf{X}=[\mathbf{X}_1, \dots, \mathbf{X}_{nr}]$ where the diagonal elements of $\mbf{R}$ are chosen to be positive real numbers. From the $\FF$-unitary invariance of $\mathbf{X}$, it follows that $\mathbf{Q}$ and $\mathbf{R}$ are independent and $\mathbf{Q}$ is a Haar distributed $\mathbb{F}$-unitary matrix. (See, e.g.,~\cite{Meckes2019}.)
By Theorem~\ref{thm:bartlett}, we also know that the diagonal entries $R_{i,i}$, $i\in \lst{rn}$, of $\mathbf{R}$ are independent with \smash{${R}_{i,i}\sim \frac{1}{\sqrt{\beta}}\chi_{\beta(rn+1-i)}$}.
From $\mbf{X}=\mbf{Q}\mbf{R}$, we have
\[
\mbf{X}_i=R_{i,i} \mbf{Q}_i+\sum_{j=1}^{i-1} R_{j,i} \mbf{Q}_j,
\]
and using elementary row operations, we get
\begin{align*}\label{eq:magic_dist}
 \det\mbf{M}({\lambda}, \mbf{X})=\prod_{a=1}^r R_{a,a}^n
\cdot \det\mbf{M}({\lambda}, \mbf{Q}).
\end{align*}
Using the independence of $\mbf{R}$ and $\mbf{Q}$ the statement of the lemma now follows.
\end{proof}

Our second observation is that, by regrouping terms in its Laplace expansion, the determinant of $\mbf{M}(\lambda, \mbf{Q})$ can be
written as a weighted sum of Vandermonde determinant products. We will index the permutations in this expansion in a particular way, and this appears frequently enough that we isolate its definition.

\begin{Definition}
 For $\underline{\cA}=(\cA_1,\dots,\cA_r)\in \cP_{r,n}$, let $\sigma_{\underline{\cA}}\in S_{rn}$
 denote the following permutation of the set $\lst{rn}$
$
\sigma_{\underline{\cA}}(a+(m-1)n)= $ the $a$-th largest element of $ \cA_m$, for $a\in \lst{n}$, $ m\in \lst{r}$.
Now when $r=2$, the elements of $(\cA_1, \cA_2)$ of $\cP_{2,n}$ can be identified just with the subsets $\cA\subset \lst{2n}$ of cardinality~$n$. Using $\cA'$ for the complement of $\cA$ in $\lst{2n}$ throughout, we further denote
\begin{equation}
\label{shorthand}
\sigma_{{\cA}} := \sigma_{\cA, \cA'} \in S_{2n}.
\end{equation}

\end{Definition}

With this set, we have the following representation.

\begin{Lemma}\label{prop:M_expansion}
 Let $\lambda_i, i \in \lst{rn}$, and $Q_{a,b}, a\in \lst{r}, b\in \lst{rn}$ be fixed numbers. For $m\in \lst{r}$, $\cA\subset \lst{rn}$, let $Q_m(\cA)=\prod_{a\in \cA} Q_{m,a}$.
 %For $\underline{\cA}=(\cA_1,\dots,\cA_r)\in \cP_{r,n}$ let $\sigma_{\underline{\cA}}$ be the permutation in $S_{rn}$ that for $m\in \lst{r}$, $a\in \lst{n}$ takes $a+(m-1)r$ into the $a$th largest element of $\cA_m$.
 Then
\begin{align*}\label{eq:detM}
 \det \mbf{M}({\lambda}, \mbf{Q})&= \sum_{(\cA_1,\dots,\cA_r)\in \cP_{r,n}} \!\!\sgn(\sigma_{\underline{\cA}}) \prod_{m=1}^r Q_{m}(\cA_m) \Delta(\cA_m).
\end{align*}
\end{Lemma}

\begin{proof}
First, any permutation $\sigma\in S_{rn}$ of $\lst{rn}$ can be mapped in a~one-to-one way to a~pair consisting of an element of $\cP_{r,n}$ and a~vector $(\sigma_1, \dots, \sigma_r)$ of $r$ permutations of $S_n$. For $\sigma\in S_{rn}$, $m\in \lst{r}$, let
$
 \cA_m=\{\sigma(i) \mid i\in \lst{rn}, \, \nfr{i}{r}=m\} =\{\sigma(m), \sigma(m+r), \dots, \sigma(m+r(n-1))\}$,
and let~$\sigma_m$ be the permutation taking the ordered version of~$\cA_m$ to the sequence
\[(\sigma(m), \sigma(m+r), \dots, \sigma(m+r(n-1))).
\]
(Strictly speaking, $\sigma_m$ is the permutation that is generated by this mapping on the relative ranks of the elements of $\cA_m$ within the set.) Note that $\sigma$ can be reconstructed from the pair~$(\cA_1,\dots, \cA_r)$, $(\sigma_1, \dots, \sigma_r)$, and we have that
$
 \sgn(\sigma)=\sgn(\sigma_{\underline{\cA}}) \prod_{m=1}^r \sgn(\sigma_m)$.
Next, recall the expression of the entries $M_{i,j}$ of $\mbf{M}=\mbf{M}({\lambda}, \mbf{Q})$ from \eqref{eq:MagicM_ij},
 and write out the determinant expansion
 \begin{align*}
 \det \mbf{M}&= \det \mbf{M}^{\dagger}=\sum_{\sigma\in S_{rn}} \sgn(\sigma) \prod_{i=1}^{rn} M_{\sigma(j),i}=\sum_{\sigma\in S_{rn}} \sgn(\sigma) \prod_{m=1}^{r} \prod_{a=0}^{n-1} \lambda_{\sigma(m+ar)}^{a} Q_{m,\sigma(m+ar)}
 \\
 &=\sum_{(\cA_1,\dots,\cA_r)\in \cP_{r,n}}\!\!\sgn(\sigma_{\underline{\cA}}) \prod_{m=1}^r Q_{m}(\cA_m) \sum_{\sigma_1, \dots, \sigma_r} \prod_{m=1}^r \sgn(\sigma_m) \prod_{a=0}^{n-1} \lambda_{\sigma_m(m+a r)}^{a}.
\end{align*}
Note that the second summation in the last line is over the possible vector of permutations $(\sigma_1,\dots,\sigma_r)\in S_n^r$ corresponding to a particular partition $(\cA_1,\dots,\cA_r)\in \cP_{r,n}$. For a fixed $m\in \lst{r}$ and $\cA_m$, the sum $\sum_{\sigma_m} \sgn(\sigma_m) \prod_{a=0}^{n-1} \lambda_{\sigma_m(m+a r)}^{a}$ is just the Vandermonde determinant~$\Delta(\cA_m)$ which
concludes the proof.
\end{proof}

\section{Evaluating the second moment of the determinant}\label{sec:2moment}

The next proposition gives the proof of the first statement of Theorem~\ref{thm:moment}.

\begin{Proposition}\label{prop:det_l=2}
Fix $r\ge 2$, $n\ge 2$. Let $\mbf{Q}$ be an $rn\times rn$ Haar distributed $\mathbb{F}$-unitary matrix.
Then, for all ${\lambda}\in \R^{rn}$ we have
 \begin{align*}%\label{eq:det2_general}
 E_{\boldsymbol{Q}} \big[|{\det \mbf{M}(
 {\lambda}, \mbf{Q})}|^2\big] = \kappa_{n,\beta, r,2} \sum_{(\cA_1,\dots, \cA_r)\in \cP_{r,n}} \prod_{j=1}^r \Delta
 (\cA_j)^2.
 \end{align*}
Here $\kappa_{n,\beta,r,2}$ is as defined in Lemma~{\rm\ref{lem:Gauss_red}}.
\end{Proposition}

We provide two proofs. The first relies on the Gaussian reduction of Lemma~\ref{lem:Gauss_red}. The second is provided for historical interest: it only uses the exchangeability properties of the Haar distribution, together with a determinantal identity from 1851 due to Sylvester~\cite{Sylvester1851}.
As the notation involved in the latter approach gets fairly heavy, we only present this second proof for~$r=2$, $\beta =1$. We also note that in the $r=2$ case one can use the identity~\eqref{eq:Pfaffian_detM} to obtain yet another independent proof.

\begin{proof}[Proof of Proposition~\ref{prop:det_l=2} using Gaussian reduction]
According to Lemma~\ref{lem:Gauss_red}, the statement is equivalent to the following.
 Let $\mbf{X}$ be an $r\times (rn)$ matrix with i.i.d.~standard real or complex Gaussian entries $X_{a,b}$, $ a\in\lst{r}$, $n\in \lst{rn}$, then
\begin{align}\label{eq:detM_2_Gauss}
 E \big[|{\det \mbf{M}({\lambda}, \mbf{X})}|^2\big] = \sum_{(\cA_1,\dots, \cA_r)\in \cP_{r,n}} \prod_{j=1}^r \Delta(\cA_j)^2.
 \end{align}
For a subset $\cA\subset \lst{rn}$ and $m\in \lst{r}$, we denote by $X_{m}(\cA)$ the product $\prod_{j\in \cA} X_{m,j}$.

From Proposition~\ref{prop:M_expansion}, we get
\begin{align*}
 E\big[|{\det \mbf{M}(\lambda, \mbf{X})}|^2\big]={}&\sum_{\underline{\cA}\in \cP_{r,n}}\sum_{\underline{\cB}\in \cP_{r,n}}\sgn(\sigma_{\underline{\cA}})\sgn(\sigma_{\underline{\cB}}) \\
 & \times E\left[\prod_{m=1}^r X_{m}(\cA_m) \overline{X}_{m}(\cB_m)\right]
\prod_{m=1}^r \Delta(\cA_m) \Delta(\cB_m).
\end{align*}
Since $X_{a,b}$ are i.i.d.~$\FF N(0,1)$ random variables, we have
\begin{align*}
 E\left[\prod_{m=1}^r X_{m}(\cA_m) \overline{X}_{m}(\cB_m)\right] =\begin{cases}
1 &\text{if $\cA_m=\cB_m$ for all $m\in \lst{r}$,}\\
0 &\text{otherwise.}
 \end{cases}
\end{align*}
From this, the identity \eqref{eq:detM_2_Gauss} and the statement of the proposition follows.
\end{proof}

\begin{proof}[Proof of Proposition~\ref{prop:det_l=2} for $\boldsymbol{ r\!=2}$ and $\boldsymbol{\beta\!=1}$ using exchangeability]
Using again Lem\-ma~\ref{prop:M_expansion}, we have that
\begin{align*}
 E\big[\det \mbf{M}(\lambda,\mbf{Q})^2\big]={}&\sum_{|\cA|=n,\, \cA\subset\lst{2n}}\sum_{|\cB|=n,\, \cB\subset\lst{2n}}\!\!
 \sgn(\sigma_{\cA})
 \sgn(\sigma_{\cB})
 E[ Q_{1}(\cA) Q_{2}(\cA')Q_{1}(\cB) Q_{2}(\cB')] \\
 & \times
 \Delta(\cA) \Delta(\cA')\Delta(\cB) \Delta(\cB').
\end{align*}
Recall the shorthand \eqref{shorthand}.
The expression $Q_{1}(\cA) Q_{2}(\cA')Q_{1}(\cB) Q_{2}(\cB')$ is a monomial of the form \smash{$\prod_{j=1}^{2n} q_{1,j}^{a_j} q_{2,j}^{2-a_j}$}, where for each $j$ we have $a_j\in \{0,1,2\}$. By the exchangeability of the columns of a Haar orthogonal matrix, we have
\begin{align}\label{eq:def_g}
 E\left[\prod_{j=1}^{2n} q_{1,j}^{a_j} q_{2,j}^{2-a_j}\right]=g(c_0,c_1,c_2)
\end{align}
for some function $g$ where $c_i=|\{a_j=i\mid j\in \lst{2n}\}|$. Moreover,
\begin{align*}
 E[Q_{1}(\cA) Q_{2}(\cA')Q_{1}(\cB) Q_{2}(\cB')]=g(a,2n-2a,a), \qquad a=|\cA\cap \cB|=|\cA' \cap \cB'|.
\end{align*}
Thus
\begin{align*}
 E\big[ |{\det \mbf{M}(\lambda,\mbf{Q})}|^2\big]={}&\sum_{|\cA|=n,\, \cA\subset\lst{2n}}\sgn(\sigma_{\cA})
 \Delta(\cA) \Delta(\cA')\\
 & \times
 \sum_{a=0}^n g(a,2n-2a, a) \sum_{|\cB|=n,\, |\cA\cap \cB|=a}
 \sgn(\sigma_{\cB})
 \Delta(\cB) \Delta(\cB').
\end{align*}
To finish, we will show that for any fixed $0\le a\le n$, $\cA\subset \lst{2n}$ with~${|\cA|=n}$, and $\cS\subset \cA$ with $|\cS|=a$ we have
\begin{align}\label{eq:det_sylv}
(-1)^{n-a} \Delta(\cA) \Delta(\cA') = \sgn(\sigma_{\cA}) \sum_{|\cB|=n,\, \cA\cap \cB=\cS}
 \sgn(\sigma_{\cB})
 \Delta(\cB) \Delta(\cB').
\end{align}
This will imply
\begin{align*}
 E\big[|{\det \mbf{M}(\lambda, \mbf{Q})}|^2\big]=\left(\sum_{a=0}^{n}(-1)^{n-a} \binom{n}{a}g(a,2n-2a,a)\right) \sum_{|\cA|=n,\, \cA\subset\lst{2n}} \Delta(\cA)^2 \Delta(\cA')^2,
\end{align*}
and so proves the proposition for $r=2$ with the constant left in the form
\begin{align*}
c_{n,2,2}&=\sum_{a=0}^{n}(-1)^{n-a} \binom{n}{a}g(a,2n-2a,a)=E\left[\prod_{j=1}^n \bigl(q_{1,2j-1}q_{2,2j-1}q_{1,2j}q_{2,2j}-q_{1,2j-1}^2 q_{2,2j}^2\bigr)\right],
\end{align*}
recall \eqref{eq:def_g}.

Now, the Vandermonde identity \eqref{eq:det_sylv} has in fact previously appeared in a
paper of Gunson~\cite{Gunson}. Here we highlight that it is also a consequence of an earlier, more general determinantal result.

First, when $a=n$ the sum on the right side
of \eqref{eq:det_sylv} only contains the term corresponding to $\cB=\cA$, so the two sides are equal. When $a=0$, then again the right side only contains one term, corresponding to $\cB=\cA'$. Since $\sigma_{\cA'}$ can be obtained from $\sigma_{\cA}$ using $n$ transpositions, \eqref{eq:det_sylv} holds in this case as well.
For the $1\le a\le n-1$ case we use the following fact due to Sylvester~\cite{Sylvester1851}.
Suppose that $\mbf{R}_1$, $\mbf{R}_2$ are $n\times n$ matrices and $1\le m\le n$. Let $I$ be a fixed subset of column indices of $\mbf{R}_1$ with $|I|=m$.
Then
\begin{align}\label{eq:Sylvester}
 \det(\mbf{R}_1) \det(\mbf{R}_2)=\sum_{(\widehat{\mbf{R}}_1, \widehat{\mbf{R}}_2)}
 \det\bigl(\widehat{\mbf{R}}_1\bigr) \det\bigl( \widehat{\mbf{R}}_2\bigr),
\end{align}
where the sum is over all pairs of $n\times n$ matrices $\bigl(\widehat{\mbf{R}}_1, \widehat{\mbf{R}}_2\bigr)$
which can be obtained from $(\mbf{R}_1, \mbf{R}_2)$ by interchanging the $m$ columns of $\mbf{R}_1$ corresponding to $I$ with $m$ columns of $\mbf{R}_2$, while preserving the ordering of the columns.

Again, fix $1\le a\le n-1$, $\cA\subset \lst{2n}$ with $|\cA|=n$, and $\cS\subset \cA$ with $|\cS|=a$. Let $\mbf{R}_1$ be the~${n\times n}$ Vandermonde matrix corresponding to $\cA$, $\mbf{R}_2$ the Vandermonde matrix corresponding to $\cA'$, and consider Sylvester's identity applied for~$m=n-a$ where the fixed set of coordinates of $\mbf{R}_1$ correspond to $\widetilde{\cS}:=\cA\setminus \cS$.
For each $\bigl(\widehat{\mbf{R}}_1,\widehat{\mbf{R}}_2\bigr)$ in the sum, let $\cB$ be the set of indices corresponding to the columns of \smash{$\widehat{\mbf{R}}_1$}, then $|\cB|=n$, $\cA\cap \cB=\cS$. We set $\cT=\cA' \cap \cB'$ and~$\smash{\widetilde{\cT}=\cA'\setminus \cT}$.
The left side of \eqref{eq:Sylvester} is $\Delta(\cA)\Delta(\cA')$, while
the sum on the right side is
\begin{align*}
 \sum_{|\cB|=n, \cA\cap \cB=\cS}\sgn(\sigma_{\cA,\widetilde \cS\to \widetilde \cT})
\sgn(\sigma_{\cA',\widetilde \cT\to \widetilde \cS})
\Delta(\cB)\Delta(\cB').
\end{align*}
Here \smash{$\sigma_{\cA,\widetilde \cS\to \widetilde \cT}$} is the permutation of the columns of $\cA$ that we get after replacing the columns corresponding to $\wt \cS$ with the columns corresponding to $\wt \cT$ from $\cA' $ (keeping their order), and then reordering them based on their indices. The permutation \smash{$\sigma_{\cA',\widetilde \cT\to \widetilde \cS}$} is defined similarly. The identity \eqref{eq:det_sylv} now comes down to showing that
\begin{align}\label{eq:perm_sgn}
\sgn(\sigma_{\cA,\widetilde \cS\to \widetilde \cT})
\sgn(\sigma_{\cA',\widetilde \cT\to \widetilde \cS})=(-1)^{n-a} \sgn(\sigma_{\cA})
 \sgn(\sigma_{\cB}).
\end{align}
Consider the permutation $\sigma_{\cA}$. Interchange the indices in $\wt \cS$ with $\wt \cT$ (keeping their respective order), since $\big|\wt \cS\big|=\big|\wt \cT\big|=n-a$, we can achieve this with $n-a$ transpositions. Then reorder the elements in positions $1,3,\dots, 2n-1$, and then in the positions $2,4, \dots, 2n$. We can do this using the permutations $\sigma_{\cA,\widetilde \cS\to \widetilde \cT}$ and $\sigma_{\cA',\widetilde \cT\to \widetilde \cS}$ applied to these collections of indices, so the signature of this step is exactly
\smash{$\sgn(\sigma_{\cA,\widetilde \cS\to \widetilde \cT})
\sgn(\sigma_{\cA',\widetilde \cT\to \widetilde \cS})$}. The resulting permutation is exactly~$\sigma_{\cB}$, and the signature of the steps we have taken is \smash{$(-1)^{n-a}\sgn(\sigma_{\cA,\widetilde \cS\to \widetilde \cT})
\sgn(\sigma_{\cA',\widetilde \cT\to \widetilde \cS})$}, which proves~\eqref{eq:perm_sgn}.
\end{proof}

\section[Evaluating fourth moment of the determinant for r=2]{Evaluating fourth moment of the determinant for $\boldsymbol{ r=2}$}
\label{sec:4moment}

 We show the following.

\begin{Proposition}\label{prop:det_l=4}
Fix $n\ge 2$. Let $\mbf{X}$ be a $2\times (2n)$ matrix of independent $\FF N(0,1)$ random variables.
% as the first $2$ rows of a size $2n$ Haar orthogonal matrix.
For all ${\lambda}\in \R^{2n}$ it holds that
 \[%\label{eq:det4}
 E |{\det \mbf{M}({\lambda}, \mbf{X})}|^4 = (12/\beta)^n \sum_{\cA\subset \lst{2n},\, |\cA|=n} \Delta (\cA)^4 \Delta(\cA')^4.
 \]
\end{Proposition}

This establishes the second statement of Theorem~\ref{thm:moment} via the Gaussian reduction Lemma~\ref{lem:Gauss_red}.

\begin{proof}
%[Proof of Proposition~\ref{prop:det_l=4}]
We carry out the proof in the $\beta=1$ case, when
the entries of $\mbf{X}$ are real standard Gaussians, commenting on the necessary modifications when $\beta =2$ at the end.

As both sides of the claimed identity are continuous in $\lambda_{i}$, $i\in \lst{2n}$, we can assume
that those variables are distinct. Next,
by Proposition~\ref{prop:M_expansion}, we have
\[
\det \mbf{M}({\lambda},\mbf{X})=\sum_{\cA\subset \lst{2n},\, |\cA|=n} \sgn(\sigma_{\cA}) X_1(\cA)X_2(\cA') \Delta(\cA)\Delta(\cA')
\]
and so also: with now $\mbf{M} = \mbf{M}(\lambda, \mbf{X})$ throughout,
\begin{align}\label{eq:det^4}
 E\big[ \det \mbf{M}^4\big]&=\sum_{\substack{\cA_i\subset \lst{2n}, \, |\cA_i|=n\\ 1\le i\le 4}} \prod_{i=1}^4 \sgn(\sigma_{\cA_i}) E \left [\prod_{i=1}^4 X_1(\cA_i) X_2(\cA_i') \right] \prod_{i=1}^4 \Delta(\cA_i)\Delta(\cA_i').
\end{align}
For a given choice of $\cA_1$, $\cA_2$, $\cA_3$, $\cA_4$ and $(j_1,j_2,j_3,j_4)\in \{0,1\}^4$, let
\begin{equation}
\cB_{(j_1,j_2,j_3,j_4)}=\cA_1^* \cap \cA_2^* \cap\cA_3^*\cap \cA_4^*, \qquad \text{with} \ \cA_i^*=\begin{cases}
 \cA_i& \text{if } j_i=1,\\
 \cA_i'& \text{if } j_i=0. \label{overlaps}
\end{cases}
\end{equation}
Also set \smash{$b_{(j_1,j_2,j_3,j_4)}=|\cB_{(j_1,j_2,j_3,j_4)}|$}.
Because the entries of $\mbf{X}$ are i.i.d.~standard normals, the expected value of $\prod_{i=1}^4 X_1(\cA_i) X_2(\cA_i')$ is zero if any entry of $\mbf{X}$ appears an odd number of times. Since the fourth moment of a standard normal is 3, this yields
\begin{align}\label{eq:det4_111}
 E\left[\prod_{i=1}^4 X_1(\cA_i) X_2(\cA_i')\right]=3^{b_{(1,1,1,1)}+b_{(0,0,0,0)}} \ind\left(b_{(j_1,j_2,j_3,j_4)}=0 \text{ if $\sum_{i=1}^4 j_i$ is odd}\right).
\end{align}
In other words, after we take expected value in \eqref{eq:det^4}, only those quadruples $(\cA_1,\dots, \cA_4)$ contribute for which intersections of the form $\cA_{i_1}'\cap \cA_{i_2}\cap \cA_{i_3} \cap\cA_{i_4}$ and $\cA_{i_1}\cap \cA_{i_2}' \cap\cA_{i_3}' \cap\cA_{i_4}'$ are empty (for distinct $i_1,\dots, i_4$).
Let $\cG_n$ denote the set of quadruples $(\cA_1,\dots, \cA_4)$ that satisfy this compatibility condition together with $|\cA_i|=n$, $\cA_i\subset \lst{2n}$.

To simplify notation, we encode the sequences $(j_1,j_2,j_3,j_4)$ where $\sum_{i=1}^4 j_i$ is even with the integers 1, \dots, 8 as follows:
\begin{align}
\label{setlist}
 \begin{array}{c|c|c|c|c|c|c|c}
 1111 & 1100 & 1010 & 1001 & 0110 & 0101 & 0011 & 0000\\ \hline
 1 & 2 & 3 & 4 & 5 & 6 & 7 & 8\\
 \end{array}
\end{align}
Hence $\cB_{(1,0,1,0)}$ and $\cB_3$ both denote the set $\cA_1 \cap\cA_2' \cap\cA_3 \cap\cA_4'$, and $b_3=b_{(1,0,1,0)}$. At the same time,
if $(\cA_1,\dots, \cA_4)\in \cG_n$, then each $\cA_i$ and $\cA_i'$ can be written as the disjoint union of four of the sets $\cB_j$, $1\le j\le 8$. For example, $\cA_1=\cB_1\cup \cB_2\cup \cB_3\cup \cB_4$.
From $|\cA_i|=|\cA_i'|$, $1\le i\le 4$, we get four equations with the sum of four terms on both sides. For instance, $|\cA_1|=|\cA_1'|$ yields
$
b_1+b_2+b_3+b_4=b_5+b_6+b_7+b_8$,
and so on.
These four equations can be readily reduced to the simpler identities
$
 b_j=b_{9-j}$, $ 1\le j\le 4$.
Equations \eqref{eq:det^4} and \eqref{eq:det4_111} together with $b_1=b_8$ lead~to
\begin{align}\label{eq:det_4_2}
 E\left[ \det \mbf{M}^4\right]=\sum_{(\cA_1,\dots,\cA_4)\in \cG_n} 9^{ b_1}
 \prod_{i=1}^4 \sgn(\sigma_{\cA_i})
 \prod_{i=1}^4 \Delta(\cA_i)\Delta(\cA_i').
\end{align}

Next, for a particular $(\cA_1,\dots, \cA_4)\in \cG_n$ consider the following ``special'' ordering of the elements of $\lst{2n}$: $x\prec y$ if $x\in \cB_i, y\in \cB_j$ with $i<j$ or if $x<y$ with $x,y\in \cB_i$. Denote by~\smash{$\wt \Delta(\cA)$} the Vandermonde determinant corresponding to the elements $\lambda_i$, $i\in \cA$, using our special ordering, and let $\sigma_i$ be the permutation that maps $1,3,\dots, 2n-1$ and $2,4,\dots, 2n$ into the elements of $\cA_i$ and $\cA_i'$ according to the special order. These permutations can be visualized using the following table:
\begin{align*}
 \begin{array}{r|c|c}
 &1,3,\dots,2n-1&2,4,\dots,2n
 \\ \hline
\sigma_1&\cB_1, \cB_2,\cB_3,\cB_4&\cB_5, \cB_6,\cB_7,\cB_8\\
\hline
\sigma_2&\cB_1, \cB_2,\cB_5,\cB_7&\cB_3, \cB_4,\cB_7,\cB_8\\
\hline
\sigma_3&\cB_1, \cB_3,\cB_5,\cB_7&\cB_2, \cB_4,\cB_6,\cB_8\\
\hline
\sigma_4&\cB_1, \cB_4,\cB_6,\cB_7&\cB_2, \cB_3,\cB_5,\cB_8
 \end{array}
\end{align*}
Note that for each $1\le i\le 4$ the permutation $\sigma_{\cA_i}$ can be obtained from $\sigma_{i}$ by reordering the elements $\sigma(2i-1)$, $1\le i\le n$, and $\sigma(2i)$, $1\le i\le n$, respectively. These two permutations correspond to the column permutations that take the Vandermonde matrices corresponding to~$\wt \Delta(\cA_i)$ and $\wt \Delta(\cA_i')$ into the Vandermonde matrices of $\Delta(\cA_i)$ and $\Delta(\cA_i')$. From this observation we get
$
 \sgn(\sigma_{\cA_i})
 \Delta(\cA_i)\Delta(\cA_i')= \sgn(\sigma_i)
 \wt \Delta(\cA_i)\wt \Delta(\cA_i')$, $1\le i\le 4$,
and hence
\begin{align*}\label{eq:det_4_3}
 E\big[ \det \mbf{M}^4\big]=\sum_{(\cA_1,\dots,\cA_4)\in \cG_n} 9^{ b_1}
 \prod_{i=1}^4 \sgn(\sigma_i)
 \prod_{i=1}^4 \wt\Delta(\cA_i)\wt\Delta(\cA_i'),
\end{align*}
as an equivalent form of \eqref{eq:det_4_2}.

One reason for introducing the above ordering is that:
for every $(\cA_1,\dots,\cA_4)\in \cG_n$,
\[
\prod_{i=1}^4 \sgn(\sigma_i)=1.
\]
To see this, first note that since $|\cB_3|=|\cB_6|=b_3$ and $|\cB_4|=|\cB_5|=b_4$, we can obtain $\sigma_2$ from $\sigma_1$ using the following steps: switch the respective elements of $\cB_3$ with $\cB_6$, and then the elements of $\cB_4$ with those of $\cB_5$ (this takes $b_3+b_4$ transpositions). Then switch the order of the `blocks'~$\cB_6$ and $\cB_5$, and the blocks~$\cB_4$ and $\cB_3$ (these two steps can be done with the same number of transpositions). We can visualize these steps as follows
\begin{align*}
 (\cB_1, \cB_2,\cB_3,\cB_4|\cB_5, \cB_6,\cB_7,\cB_8 ) &\to (\cB_1, \cB_2,\cB_6,\cB_5|\cB_4, \cB_3,\cB_7,\cB_8 ) \\
&\to (\cB_1, \cB_2,\cB_5,\cB_6|\cB_3, \cB_4,\cB_7,\cB_8 ).
\end{align*}
This shows that
\begin{align}
 \sgn(\sigma_1)=\sgn(\sigma_2) (-1)^{b_3+b_4}.\label{eq:sign_1}
\end{align}
Similarly, we can obtain $\sigma_4$ from $\sigma_3$ by switching the respecting elements of $\cB_i$ with $\cB_{9-i}$ for~${i=3,4}$, and then switch the order of $\cB_6$, $\cB_4$ and $\cB_5$, $\cB_3$
\begin{align*}
 (\cB_1, \cB_3,\cB_5,\cB_7|\cB_2, \cB_4,\cB_6,\cB_8 )& \to (\cB_1, \cB_6,\cB_4,\cB_7|\cB_2, \cB_5,\cB_3,\cB_8 ) \\
& \to (\cB_1, \cB_4,\cB_6,\cB_7|\cB_2, \cB_3,\cB_5,\cB_8 ).
\end{align*}
This yields
\begin{align}\label{eq:sign_2}
 \sgn(\sigma_3)=\sgn(\sigma_4) (-1)^{b_3+b_4}.
 \end{align}
Combining \eqref{eq:sign_1} and \eqref{eq:sign_2} produces the claimed $\prod_{i=1}^4 \sgn(\sigma_i)=1$.

{\samepage
At this point then we have the identity,
\begin{align}
 \label{eq:det_4_4}
 E\big[ \det \mbf{M}^4\big]=\sum_{(\cA_1,\dots,\cA_4)\in \cG_n} 9^{b_1}
 \prod_{i=1}^4 \wt\Delta(\cA_i)\wt\Delta(\cA_i'),
\end{align}
and our next step is to rewrite the right-hand side of \eqref{eq:det_4_4} as a product of Cauchy determinants.}

Towards this, for two disjoint sets $\cS_1, \cS_2\subset \lst{2n}$ we introduce the notation
\begin{align}\label{eq:defdelta}
 \dd(\cS_1,\cS_2):= \prod_{s_1\in \cS_1} \prod_{s_2\in \cS_2}(\lambda_{s_1}-\lambda_{s_2}),
\end{align}
with the empty product defined as 1.
For a particular $(\cA_1,\dots, \cA_4)\in \cG_n$ for each $1\le i\le 4$, we can group the variables into the sets $\cB_j$ and write
\begin{align*}%\label{eq:prod_det_0}
\wt\Delta(\cA_i)\wt\Delta(\cA_i')=\prod_{j=1}^8 \Delta(\cB_j)\prod_{a_1>a_2} \dd(\cB_{a_1}, \cB_{a_2})\prod_{b_1>b_2} \dd(\cB_{b_1}, \cB_{b_2}).
\end{align*}
Here $a_1$, $a_2$ are indices from the ``$\cB$-partition'' of $\cA_i$ and $b_1$, $b_2$ are indices from the ``$\cB$-partition'' of $\cA_i'$. Multiplying these equations for $1\le i\le 4$
and some bookkeeping yields
\begin{align}\label{eq:prod_det_44}
 \prod_{i=1}^4\wt\Delta(\cA_i)\wt\Delta(\cA_i')=\prod_{j=1}^8 \Delta(\cB_j)^4 \prod_{\substack{9- j_1\neq j_2\\1\le j_2< j_1\le 8}} \dd(\cB_{j_1}, \cB_{j_2})^2.
\end{align}
A similar argument gives
\begin{align}\label{eq:Delta_prod}
 \Delta(\lambda)^2=\prod_{j=1}^8 \Delta(\cB_j)^2 \prod_{1\le a<b\le 8} \dd(\cB_b, \cB_a)^2.
\end{align}
Recall next the Cauchy determinant formula. For two disjoint sets $\cS_1, \cS_2\subset \lst{2n}$ with $|\cS_1|=|\cS_2|$,
we denote
\begin{align*}%\label{eq:Cauchy_det}
 C_{\cS_1,\cS_2}=\det \left(\frac{1}{\lambda_{s_1}-\lambda_{s_2}}\right)_{s_1\in \cS_1, s_2, \cS_2}
 = (-1)^{\binom{|\cS_1|}{2}} \cdot \frac{\Delta(\cS_1)\Delta(\cS_2)}{\dd(\cS_1,\cS_2)}.
\end{align*}
(If $\cS_1$, $\cS_2$ are empty, we set $C_{\cS_1,\cS_2}=1$.)
From \eqref{eq:prod_det_44} and \eqref{eq:Delta_prod}, we
may deduce that
\begin{align}\label{eq:prod_det_4}
 \Delta(\lambda)^{-2} \prod_{i=1}^4\wt\Delta(\cA_i)\wt\Delta(\cA_i')=\prod_{j=1}^4 C_{\cB_j,\cB_{9-j}}^2,
\end{align}
again for any $(\cA_1,\dots, \cA_4)\in \cG_n$

The point is that by Lemma~\ref{lem:Cauchy_cycle} below we have the following. Suppose that $b_i$, $1\le i\le 8$ are fixed nonnegative integers with $b_i=b_{9-i}$ and $\sum_{i=1}^4 b_i=n$. Let $\cP_{\underline{b}}$ be the collection of partitions~${\underline{\cB}=(\cB_1,\dots, \cB_8)}$ of $\lst{2n}$ with $|\cB_i|=b_i$, $1\le i\le 8$. Then
\begin{align}\label{eq:Cauchy_sum}
 \sum_{\underline{\cB}\in \cP_{\underline{b}}} \prod_{i=1}^4 C_{\cB_i,\cB_{9-i}}^2 =\binom{n}{b_1,b_2,b_3,b_4} \sum_{\cA\subset \lst{2n}, |\cA|=n} C_{\cA,\cA'}^2.
\end{align}
Taking \eqref{eq:Cauchy_sum} for granted, \eqref{eq:prod_det_4}
allows us to rewrite \eqref{eq:det_4_4}
as in
\begin{align}\notag
 E\big[ \det \mbf{M}^4\big]&=\Delta(\lambda)^2\sum_{b_1+\cdots+b_4=n,\, b_i\ge 0} 9^{b_1} \binom{n}{b_1,b_2,b_3,b_4} \sum_{\cA\subset \lst{2n},\, |\cA|=n} C_{\cA,\cA'}^2\\
 &=12^n \times \Delta(\lambda)^2 \sum_{\cA\subset \lst{2n},\, |\cA|=n} C_{\cA,\cA'}^2\label{eq:det4_9},
\end{align}
having used the multinomial theorem in the second line.

We can now finish the proof. For any $\cA\subset \lst{2n}$ with $|\cA|=n$,
$%\label{eq:det2n_C}
\Delta(\lambda)^2\!=\! \Delta(\cA)^2 \Delta(\cA')^2 \dd(\cA,\cA')^2$,
which by the Cauchy determinant formula implies
$
\Delta(\lambda)^2 C_{\cA,\cA'}^2= \Delta(\cA)^4 \Delta(\cA')^4$.
%\dd(\cA,\cA')^4.\]
Substituting the above into \eqref{eq:det4_9} produces
\[
 E\left[ \det \mbf{M}^4\right]=12^n \sum_{\cA\subset \lst{2n},\, |\cA|=n} \Delta(\cA)^4 \Delta(\cA')^4,
\]
which is the claimed statement for $\beta=1$.

For $\beta =2$, the key change comes in equation \eqref{eq:det4_111}, where one now gets
\begin{gather}
E \left [\prod_{i=1}^2 X_1(\cA_i) X_2(\cA_i') \prod_{i=3}^4 \overline{X}_1(\cA_i) \overline{X}_2(\cA_i')\right]\nonumber\\
\qquad \notag=2^{b_{(1,1,1,1)}+b_{(0,0,0,0)}} \ind(b_{(j_1,j_2,j_3,j_4)}=0 \text{ if $j_1+j_2\neq j_3+j_4$}).\label{eq:new_det4}
\end{gather}
The coefficient $2$ and the updated condition in the indicator function result from the fourth absolute moment and the rotational invariance of the standard complex Gaussian distribution, respectively. Equation \eqref{eq:new_det4} now implies that if the sets $\cB_{(1,0,1,0)}$ and
$\cB_{(0,1,0,1)}$ are not empty then the corresponding quadruple $(\cA_1,\dots,\cA_4)$ does not contribute in the sum analogous to~\eqref{eq:det^4}. This translates to having $b_3=b_6=0$ in all subsequent calculations. In particular,
 the multinomial factor in \eqref{eq:det4_9} is replaced by $ 4^{b_1} \binom{n}{b_1,b_2,b_4}$, which
explains the reported constant of~$6^n$ for~${\beta=2}$.
\end{proof}

\section{Determinantal identities}
\label{sec:det_identities}

The main goal of the section is to prove the following identity which is crucial in the proof of Proposition~\ref{prop:det_l=4}.

\begin{Lemma}\label{lem:Cauchy_cycle}
Let $n\ge 2$, and suppose that $b_i$, $1\le i\le 8$, are fixed nonnegative integers with $b_i=b_{9-i}$ and $\sum_{i=1}^4 b_i=n$. Let $\cP_{\underline{b}}$ be the set of partitions $\underline{\cB}=(\cB_1,\dots, \cB_8)$ of $\lst{2n}$ with $|\cB_i|=b_i$, $1\le i\le 8$. Then
\begin{align}\label{eq:Cauchy_sum_1}
 \sum_{\underline{\cB}\in \cP_{\underline{b}}} \prod_{j=1}^4 C_{\cB_j,\cB_{9-j}}^2 =\binom{n}{b_1,b_2,b_3,b_4} \sum_{\cA\subset \lst{2n},\, |\cA|=n} C_{\cA,\cA'}^2.
\end{align}
\end{Lemma}

We also establish the following which unifies the various descriptions of our $r=2$ eigenvalue densities.

\begin{Lemma}\label{lem:det4_identities}
 Fix $n\ge 2$, and let ${\lambda}\in \R^{2n}$. Then the following identities hold:
 \begin{align*}%\label{eq:det4_id_1}
 \sum_{\cA\subset \lst{2n}, \, |\cA|=n} \Delta (\cA)^4 \Delta (\cA')^4
 & =2^{-n} \left(\sum_{\cA\subset \lst{2n},\, |\cA|=n} \Delta (\cA)^2 \Delta (\cA')^2\right)^2 \\
 & =2^n \Delta(\lambda)^2 \left|\Pf\left(\frac{\ind_{i\neq j}}{\lambda_i-\lambda_j}\right)\right|^2. % \label{eq:det4_id_2}
 \end{align*}
\end{Lemma}

First, though we present a conjectured generalization of Lemma~\ref{lem:Cauchy_cycle} motivated by the proof of Proposition~\ref{prop:det_l=4}.

\subsection{A conjecture on certain Vandermonde sums}

The Gaussian reduction was key to computing $E_{\boldsymbol{Q}} \bigl[|{\det \mbf{M}({\lambda}, \mbf{Q})}|^4\bigr]$
as it vastly cuts down the number of terms on the right-hand side of the basic formula \eqref{eq:det^4}. Assuming however that the end result is only a function of the exchangeability of the
rows and columns of $\mbf{Q}$,
leads to an interesting set of conjectures involving sums of squared Vandermonde determinants.

To explain, recall the enumeration of ``overlaps'' $\cB_{(j_1,j_2, j_3, j_4)}$ defined in \eqref{overlaps}. Here $j_i \in \{0,1\}$, and previously all overlaps for which \smash{$\sum_{i=1}^4 j_i$} is odd immediately dropped from consideration. Restoring these, we have sixteen possible overlaps $(j_1, j_2, j_3, j_4)$ which we again
list in reverse lexicographic order, and now denote their sizes by $x_i$, $i \in \{1, \dots, 16\}$. (Compare with \eqref{setlist}.) In particular, we now have, for example,
$\cB_{(1,1,0,1)}$
as the third element on the list for which we set~${x_3 = |B_{(1,1,0,1)}| }$. With these conventions, we can state the following.

\begin{Conjecture}
Fix nonnegative integers $x_1, \dots, x_{16}$ with $x_1+\dots+x_8=n$.
Then
% For any collection $S$ of fixed overlap sizes $x_1, \dots, x_{16}$ defined above:
\begin{align}
 &\sum_{(\cA_1,\dots,\cA_4)\subset \mathcal{S}_{\mbf{x}}} \prod_{i=1}^4 \sgn(\sigma_{\cA_i}) \Delta_{\cA_i}^2 \Delta_{\cA_i'}^2 \nonumber\\
&\qquad = \ind(x_i=x_{17-i})
 (-1)^{x_2+x_3+x_5+x_7} { n \choose x_1, \dots , x_8 }
 \sum_{\cA \subset \lst{2n}, |\cA|=n}
 \Delta_{\cA}^4 \Delta_{\cA'}^4.
\label{Conjecture}
\end{align}
Here $\mathcal{S}_{\mbf{x}}$ is the collection of all quadruples $(\cA_1,\dots,\cA_4)$, where the overlap sizes are given by $x_1,\dots, x_{16}$.
\end{Conjecture}

The identity \eqref{eq:Cauchy_sum_1} of Lemma~\ref{lem:Cauchy_cycle}, which can be recast in terms of Vandermonde instead of Cauchy determinants, is in fact the special case of \eqref{Conjecture} where one restricts the sum on the left-hand side to ``even" symmetric overlaps (for which $x_2 = x_3 = x_5=x_8 = 0$).

\subsection[Proofs of Lemmas 6.1 and 6.2]{Proofs of Lemmas \ref{lem:Cauchy_cycle} and
\ref{lem:det4_identities}}

We start by recalling
 the definitions of the Pfaffian and Hafnian of a matrix.

\begin{Definition}\label{def:PfaffHaf}
For a $2n \times 2n$ skew-symmetric matrix $M$, the Pfaffian of $M$ is defined by
\begin{align}\label{eq:Pfaffiandef}
 \Pf(M)=\frac{1}{2^n n!} \sum_{\sigma\in S_{2n}} \sgn(\sigma) \prod_{i=1}^n M_{{\sigma(2i-1)},{\sigma(2i)}}
 = \sum_{\alpha\in \Pi_{2n}} \sgn(\pi_\alpha) \prod_{i=1}^n M_{\alpha_{i,1},\alpha_{i,2}},
\end{align}
while for a $2n \times 2n$ symmetric matrix $M$, the Hafnian of $M$ is defined by
\begin{align}\label{eq:Hafniandef}
 \Hf(M)=\frac{1}{2^n n!} \sum_{\sigma\in S_{2n}} \prod_{i=1}^n M_{{\sigma(2i-1)},{\sigma(2i)}}
 =\sum_{\alpha\in \Pi_{2n}} \prod_{i=1}^n M_{\alpha_{i,1},\alpha_{i,2}}.
\end{align}
Here $\Pi_{2n}$ refers to the set of matchings on $\lst{2n}$, i.e., partitions $(\alpha_{1,1},\alpha_{1,2}), \dots, (\alpha_{n,1},\alpha_{n,2})$ of $\lst{2n}$ with $\alpha_{i,1}<\alpha_{i,2}$ and $\alpha_{1,1}<\cdots<\alpha_{n,1}$. In the Pfaffian formula, $\pi_\alpha$ denotes the corresponding permutation $(\alpha_{1,1},\alpha_{1,2}, \alpha_{2,1},\alpha_{2,2},\dots, \alpha_{n,1},\alpha_{n,2})$.
\end{Definition}

Breaking the proof of Lemma~\ref{lem:det4_identities} into two parts, we start with the following.

\begin{Proposition}\label{prop:det2+Pf}
Let ${\lambda}\in \R^{2n}$. Then we have
 \begin{align*}%\label{eq:det2+Pf}
 \sum_{|\cA|=n,\, \cA\subset\lst{2n}} \Delta(\cA)^2 \Delta (\cA')^2&=(-2)^n \Delta (\lambda) \Pf\left(\frac{\ind_{i\neq j}}{\lambda_i-\lambda_j}\right).
 \end{align*}
The right side is defined as the appropriate limit if the $\lambda_i$'s are not all distinct.
\end{Proposition}

\begin{proof} Assuming that the $\lambda_i$'s are distinct, the definition \eqref{eq:Pfaffiandef} gives
\begin{align*}
 \Pf\left(\frac{\ind_{i\neq j}}{\lambda_i-\lambda_j}\right)&
 %=\frac{1}{2^n n!} \sum_{\sigma\in S_{2n}} \sgn(\sigma) \prod_{i=1}^n
 %\frac{1}{\lambda_{\sigma(2i-1)}-\lambda_{\sigma(2i)}}
 =\sum_{\alpha\in \Pi_{2n}} \sgn(\pi_\alpha) \prod_{i=1}^n \frac{1}{\lambda_{\alpha_{i,1}}-\lambda_{\alpha_{i,2}}}.
\end{align*}
A matching $\alpha \in \Pi_{2n}$ can be encoded with a pair $(\cA, \tau)$ where $|\cA|=n$, $\cA\subset \lst{2n}$ and $\tau\in S_n$. Indeed, $\cA=\{\alpha_{1,1}, \dots, \alpha_{n,1}\}$ and the permutation $\tau\in S_n$ that takes the ordered version of $\cA'$ into $(\alpha_{1,2}, \dots, \alpha_{n,2})$ identifies the matching $\alpha$. Note that $\sgn(\pi_\alpha)$ is equal to $\sgn(\sigma_{\cA} )\sgn(\eta_n)$ where $\sigma_{\cA}$ is defined in \eqref{shorthand}, and $\eta_n$ is the permutation $(1,n+1,2, n+2, \dots, n, 2n)$. This leads~to
\begin{align}
 \label{eq:Pf4_2}
\Pf\left(\frac{\ind_{i\neq j}}{\lambda_i-\lambda_j}\right)=\sgn(\eta_n) \sum_{|\cA|=n,\, \cA\subset \lst{2n}} \sgn(\sigma_{\cA}) \sum_{\tau\in S_n} \sgn(\tau) \prod_{i=1}^n \frac{1}{\lambda_{a'_{\tau(i)}}-\lambda_{a_i}},
\end{align}
where we used $a_i$, $a_i'$, $i\in \lst{n}$, to denote the (ordered) elements of $\cA$ and $ \cA'$, respectively.

For a given $|\cA|=n$, $\cA\subset \lst{2n}$, we have
\[
\Delta(\lambda)=\Delta(\cA)\Delta(\cA') \prod_{\substack{1\le i<j\le 2n\\ \left|\{i,j\}\cap \cA\right|=1}} (\lambda_j-\lambda_i).
\]
By \eqref{eq:defdelta}, we have
\[
\prod_{\substack{1\le i<j\le 2n\\ \left|\{i,j\}\cap \cA\right|=1}} (\lambda_j-\lambda_i)=(-1)^{k_\cA}
\dd(\cA, \cA')
\]
 with $k_{\cA}$ the number of pairs $(i,j)$ with $1\le i<j\le 2n$ and $i\in \cA$, $j\in \cA'$. Observe then that~${n^2-k_{\cA}}$ is exactly the number of inversions of the permutation $\sigma_{\cA}$.
Hence,
 \begin{align}
 \label{eq:det_exp_A_Ac}
\Delta
%_{{\lambda}}
(\lambda)=(-1)^n \sgn(\sigma_{\cA}) \Delta
%_{{\lambda}}
(\cA)\Delta
%_{{\lambda}}
(\cA') \dd(\cA, \cA').
\end{align}
Using this identity, we get
\begin{align}
 \sum_{|\cA|=n,\, \cA\subset\lst{2n}} \frac{\Delta
 %_{{\lambda}}
 (\cA)^2 \Delta
 %_{{\lambda}}
 (\cA')^2}{ \Delta(\lambda)}&=(-1)^n \sum_{|\cA|=n,\, \cA\subset\lst{2n}}
 \sgn(\sigma_{\cA})
 \frac{\Delta(\cA) \Delta(\cA')}{\dd(\cA, \cA')}\nonumber \\
 &=(-1)^{n+\binom{n}{2}}\sum_{|\cA|=n,\, \cA\subset\lst{2n}} \sgn(\sigma_{\cA}) C_{\cA, \cA'}
\nonumber
 \\
 &=(-1)^{n+\binom{n}{2}}\sum_{|\cA|=n,\, \cA\subset\lst{2n}} \sgn(\sigma_{\cA}) \sum_{\tau\in S_n} \sgn(\tau) \frac{1}{\lambda_{a_i}-\lambda_{a'_{\tau(i)}}}, \label{eq:det4_123}
\end{align}
where in the last line we expanded the Cauchy determinant, and used that $a_i$, $a_i'$, $i\in \lst{n}$, for the ordered elements of $\cA$, $\cA'$, respectively.
Note that \smash{$\sgn(\eta_n)=(-1)^{\binom{n}{2}}$}, hence by comparing~\eqref{eq:det4_123} with~\eqref{eq:Pf4_2} the statement follows.
\end{proof}

\begin{Proposition}\label{prop:det4+det2}
 Let ${\lambda}\in \R^{2n}$. Then we have
 \begin{align*}%\label{eq:det4+det2}
 \sum_{\cA\subset \lst{2n},\, |\cA|=n} \Delta (\cA)^4 \Delta (\cA')^4&=2^{n} \Delta(\lambda)^2
\Hf\left(\frac{\ind_{i\neq j}}{(\lambda_i-\lambda_j)^2}\right).
 \end{align*}
 Again, the right side is defined as the appropriate limit if the $\lambda_i$'s are not all distinct.
\end{Proposition}

To prove Proposition~\ref{prop:det4+det2}, we will use the
 following lemma, which is well-known in the statistical physics community. (See, e.g., \cite[equation~(D.29)]{FGM_1995}.) We include the proof for completeness.
\begin{Lemma}[cycle cancellation]\label{lem:cycle}
Suppose that $k\ge 3$, and $z_1, \dots, z_k$ are distinct numbers. Then
\begin{align}\label{eq:cycle}
\sum^*_{\sigma} \prod_{i=1}^k \frac{1}{z_i-z_{\sigma(i)}}=0,
\end{align}
where the sum is over all cycles $\sigma$ supported on $\{1,2,\dots,k\}$. We also have
\begin{align}\label{eq:cycle_full}
 \sum_{\sigma\in S_k} \prod_{i=1}^k \frac{1}{z_{\sigma(i)}-z_{\sigma(i+1)}}=0,
\end{align}
where $\sigma(k+1)=\sigma(1)$.
\end{Lemma}
\begin{proof}
We have
\begin{align}\label{eq:cycle_3}
 \frac{1}{(z_1-z_a)(z_b-z_1)}=\frac{1}{z_b-z_a}\cdot \left(\frac{1}{z_b-z_1}-\frac{1}{z_a-z_1}\right).
\end{align}
Fix a length $k-1$ cycle $(b_1, \dots, b_{k-1})$ on $\{1,\dots, k-1\}$, and consider all length $k$ cycles that we get by inserting $k$ at some point. Then by \eqref{eq:cycle_3} the contribution of all of these cycles in the sum in \eqref{eq:cycle} (using $b_k=b_1$) is
\[
\left(\sum_{j=1}^{k-1} \frac{1}{z_{b_j}-z_k}-\frac{1}{z_{b_{j+1}}-z_k}\right) \prod_{j=1}^{k-1} \frac{1}{z_{b_j}-z_{b_{j+1}}}=0.
\]
This proves \eqref{eq:cycle}. The identity \eqref{eq:cycle_full} now also follows by observing that the sum in \eqref{eq:cycle_full} is exactly $k$ times the sum in \eqref{eq:cycle}.
\end{proof}

\begin{proof}[Proof of Proposition~\ref{prop:det4+det2}]
It is sufficient to prove the statement in the case where all the $\lambda_i$'s are distinct.

We first write
\begin{align}
\sum_{|\cA|=n,\, \cA\subset\lst{2n}} \frac{\Delta(\cA)^4 \Delta(\cA')^4}{ \Delta(\lambda)^2}& =
\sum_{|\cA|=n, \,\cA\subset\lst{2n}} C_{\cA, \cA'}^2\nonumber
\\
&=\frac{1}{(n!)^2} \sum_{\substack{a_1,\dots, a_n, a_1', \dots, a_n'\\\{a_1,\dots, a_n, a_1', \dots, a_n'\}=\lst{2n}}} C_{(a_1,\dots,a_n), (a_1', \dots, a_n')}^2. \label{eq:det4det2_0}
\end{align}
Here we used \eqref{eq:det_exp_A_Ac}, and in the second step we reordered the rows and columns of the Cauchy determinants in all possible ways.
 Next,
we expand the Cauchy determinants to get
\begin{gather}
 \sum_{|\cA|=n,\, \cA\subset\lst{2n}} \frac{
 \Delta
 %_{{\lambda}}
 (\cA)^4 \Delta
 %_{{\lambda}}
 (\cA')^4}{ \Delta
 %_{{\lambda}}
 (\lambda)^2}\nonumber\\
 \qquad=\frac{1}{(n!)^2}\sum_{a_i, a_i'} \sum_{\sigma, \tau\in S_n} \sgn\bigl(\sigma\circ \tau^{-1}\bigr)
 \prod_{i=1}^n \frac{1}{(\lambda_{a_i} - \lambda_{a'_{\sigma(i)}})(\lambda_{a_i} - \lambda_{a'_{\tau(i)}})}.\label{eq:det4_sum}
\end{gather}
(The restrictions on $a_i$, $ a_i'$ are the same as in \eqref{eq:det4det2_0}.)

We first evaluate the diagonal part of the double sum in \eqref{eq:det4_sum} by showing that
\begin{align}\label{eq:det4_diag}
\frac{1}{(n!)^2}\sum_{a_i, a_i'} \sum_{\sigma\in S_{n}}
 \prod_{i=1}^n \frac{1}{(\lambda_{a_i} - \lambda_{a'_{\sigma(i)}})^2}=2^n \Hf \left(
 \frac{\ind_{i \neq j }}{(\lambda_i -\lambda_j)^{2}} \right).
\end{align}
Consider the mapping that takes a particular choice of $a_1,\dots, a_n, a_1', \dots, a_n'$ and $\sigma\in S_n$ into the matching of $\lst{2n}$ that matches $a_i$ with \smash{$a_{\sigma(i)}'$} for $1\le i\le n$. Each particular matching $\alpha\in \Pi_{2n}$ shows up exactly $2^n (n!)^2$ times as the result of this mapping, which proves \eqref{eq:det4_diag}.

To complete the proof of our proposition, we just need to show that the `non-diagonal' terms in the sum \eqref{eq:det4_sum} cancel out. For this, it is sufficient to show that if $\sigma\neq \tau$ are fixed elements of~$S_n$ then
\begin{equation}
\label{eq:off_diag_vanish}
\sum_{\substack{a_1,\dots, a_n, a_1', \dots, a_n'\\\{a_1,\dots, a_n, a_1', \dots, a_n'\}=\lst{2n}}}
 \prod_{i=1}^n \frac{1}{(\lambda_{a_i} - \lambda_{a'_{\sigma(i)}})(\lambda_{a_i} - \lambda_{a'_{\tau(i)}})}=0.
\end{equation}
For a particular pair $\sigma$, $\tau$, consider the permutation on $\lst{2n}=\{a_1,\dots, a_n, a_1', \dots, a_n'\}$ that takes~$a_i$ to $a_{\sigma_i}'$ and~$a_i'$ to $a_{\tau^{-1}(i)}$. This permutation has even cycles, and because $\sigma\neq \tau$ it must have at least one cycle of length at least~4. Let one of these cycles be
\begin{align*}%\label{eq:cycle_123}
 a_{i_1}\to a'_{j_1} \to a_{i_2} \to a'_{j_2} \to \dots\to a'_{j_k}\to a_{i_1}.
\end{align*}
Here $1\le i_1,\dots,i_k\le n$ are distinct, and the same holds for $j_1, \dots, j_k$. Let $\cB\subset \lst{2n}$ with ${|\cB|=2k}$. We claim that
\begin{align}\label{eq:off_diag_vanish_1}
\sum_{\substack{a_1,\dots, a_n, a_1', \dots, a_n'\\\{a_1,\dots, a_n, a_1', \dots, a_n'\}=\lst{2n}\\ \{a_{i_1},a'_{j_1}, a_{i_2}, a'_{j_2} \dots, a_{i_k}, a'_{j_k}\}=\cB }}
 \prod_{i=1}^n \frac{1}{(\lambda_{a_i} - \lambda_{a'_{\sigma(i)}})(\lambda_{a_i} - \lambda_{a'_{\tau(i)}})}=0.
\end{align}
The sum on the left can be rewritten as a double sum where we first sum over all possible assignments of the values of $a_{i_1},a'_{j_1}, a_{i_2}, a'_{j_2}, \dots, a_{i_k}, a'_{j_k}$, and then in the second sum we sum over the remaining $2n-2k$ variables. Factoring out the terms corresponding to the indices not in~$\cB$, the inner sum becomes
\[
\sum_{\{a_{i_1},a'_{j_1}, a_{i_2}, a'_{j_2} \dots, a_{i_k}, a'_{j_k}\}=\cB} \prod_{\ell=1}^k \frac{1}{\lambda_{a_{i_\ell}}-\lambda_{a'_{j_\ell}}}\cdot \frac{1}{\lambda_{a'_{j_\ell}}-\lambda_{a_{i_{\ell+1}}}}
\]
with $i_{k+1}=i_1$. But this is equal to zero by \eqref{eq:cycle_full} of Lemma~\ref{lem:cycle}, as applied to the sum over the permutations of elements of $\cB$ and the values \smash{$\lambda_{a_{i_1}}, \lambda_{a'_{j_1}}, \dots, \lambda_{a_{i_k}}, \lambda_{a'_{j_k}}$}. This proves \eqref{eq:off_diag_vanish_1}, which gives \eqref{eq:off_diag_vanish} and the statement of the proposition.
\end{proof}

We now have all the components to prove Lemma~\ref{lem:det4_identities}.

\begin{proof}[Proof of Lemma~\ref{lem:det4_identities}]
The proof follows from the statements of Proposition~\ref{prop:det2+Pf} and \ref{prop:det4+det2} once we
establish that
\begin{align*}\label{eq:Pf2=Hf}
 \Pf\left(\frac{\ind_{i\neq j}}{\lambda_i-\lambda_j}\right)^2= \det\left(\frac{\ind_{i\neq j}}{\lambda_i-\lambda_j}\right) = \Hf\left(\frac{\ind_{i\neq j}}{(\lambda_i-\lambda_j)^2}\right).
\end{align*}
The first equality here is due the standard fact that
 square of the Pfaffian of a skew-symmetric matrix is equal to its determinant.
The second has in fact been noticed before in \cite{DSZ}, and can be seen by expansion
\[
 \det\left(\frac{\ind_{i\neq j}}{\lambda_i-\lambda_j}\right)=\sum_{\sigma\in S_{2n},\, \sigma(j)\neq j} \sgn(\sigma) \prod_{j=1}^{2n} \frac{1}{\lambda_j-\lambda_{\sigma(j)}},
\]
where the sum is over all the permutations of $\lst{2n}$ that have no fixed elements. Because of Lemma~\ref{lem:cycle}, the contribution of the permutations that have a cycle of length at least~3 cancels out. The remaining terms correspond to the permutations that only have 2-cycles, that is, the permutations that corresponding to perfect matchings of $\lst{2n}$. For such a permutation $\sigma$, we~have
\[
\sgn(\sigma) \prod_{j=1}^{2n} \frac{1}{\lambda_j-\lambda_{\sigma(j)}}=\prod_{i\in \cA} \frac{1}{(\lambda_{i}-\lambda_{\sigma(i)})^2},
\]
where $\cA$ is a set containing an element from each 2-cycle. By \eqref{eq:Hafniandef}, the sum of these terms is exactly the Hafnian of the $n\times n$ matrix with $(i,j)$ entry \smash{$\frac{\ind_{i\neq j}}{(\lambda_i-\lambda_j)^2}$}.
\end{proof}

We are now ready to prove Lemma~\ref{lem:Cauchy_cycle}.

\begin{proof}[Proof of Lemma~\ref{lem:Cauchy_cycle}]
 We first note that by \eqref{eq:det4det2_0} and Proposition~\ref{prop:det4+det2}, we have
 \begin{align}\label{eq:C_Hf}
 \sum_{\cA\subset \lst{2n}, |\cA|=n} C_{\cA,\cA'}^2=2^n \Hf\left(\frac{\ind_{i\neq j}}{(\lambda_i-\lambda_j)^2}\right).
 \end{align}
Fix $b_1, \dots, b_4$ with $\sum_{i=1}^4 b_i=n$ and set $b_i=b_{9-i}$ for $5\le j\le 8$. For a particular $\sigma\in S_{2n}$ denote by $\cB_{\sigma,1}, \dots, \cB_{\sigma,4},\cB'_{\sigma,1}, \dots, \cB'_{\sigma,4}$ the ordered lists we obtain once we partition $(\sigma(1),\dots, \sigma(2n))$ into parts of lengths $b_1$, $b_2$, $b_3$, $b_4$, $b_8$, $b_7$, $b_6$, $b_5$. Then we have
\begin{align*}%\label{eq:Cauchy_B4_1}
 \sum_{\underline{\cB}\in \cP_{\underline{b}}} \prod_{i=1}^4 C_{\cB_i,\cB_{9-i}}^2 =\frac{1}{(b_1! b_2! b_3! b_4!)^2}\sum_{\sigma\in S_{2n}} \prod_{j=1}^4 C_{\cB_{\sigma,j},\cB'_{\sigma,j}}^2.
\end{align*}
We introduce the temporary notation $z_i=\lambda_{\sigma(i)}$ and $z'_i=\lambda_{\sigma(n+i)}$ for $1\le i\le n$, and also $\cD_1=\{1,\dots, b_1\}$, $\cD_2=\{b_1+1,\dots, b_1+b_2\}$, $\cD_3=\{b_1+b_2+1, \dots, b_1+b_2+b_3\}$, $\cD_4=\{b_1+b_2+b_3+1,\dots,n\}$.

For a given $1\le j\le 4$, we can expand the square of the appropriate Cauchy determinant as
\begin{align*}
 C^2_{\cB_{\sigma,j},\cB'_{\sigma,j}}=\sum_{\eta_j, \tilde \eta_j\in S(\cD_j)} \sgn(\eta_j)\sgn(\tilde \eta_j) \prod_{i\in \cD_j} \frac{1}{z_i-z'_{\eta_j(i)}} \cdot \frac{1}{z_i-z'_{\tilde \eta_j(i)}},
\end{align*}
where $S(\cD_j)$ is the set of permutations of $\cD_j$.
We can represent $\eta_1,\dots, \eta_4$ as a single permutation of $\lst{n}$ that preserves $\cD_1$, $\cD_2$, $\cD_3$, $\cD_4$. Denoting the set of these permutations $S_{b_1,b_2,b_3,b_4}$, we get
\begin{align}
 \sum_{\underline{\cB}\in \cP_{\underline{b}}} \prod_{j=1}^4 C_{\cB_j,\cB_{9-j}}^2
 ={}&\frac{1}{(b_1! b_2! b_3! b_4!)^2}\sum_{\sigma\in S_{2n}} \sum_{\eta,\tilde \eta \in S_{b_1,b_2,b_3,b_4}}
 \sgn(\eta) \sgn(\tilde \eta) \nonumber\\
 & \times \prod_{i=1}^n \frac{1}{z_i-z'_{\eta(i)}}\cdot \frac{1}{z_i-z'_{\tilde \eta(i)}}.\label{eq:Cauchy_B4_1_2}
\end{align}
Just as in previous computations, we consider the diagonal and off-diagonal terms of the resulting sum separately, and show that the off-diagonal terms cancel.

The diagonal terms correspond to the cases $\eta=\tilde \eta$, this gives
\begin{align}\label{eq:Cauchy_B4_2}
 \frac{1}{(b_1! b_2! b_3! b_4!)^2}\sum_{\sigma\in S_{2n}} \sum_{\eta \in S_{b_1,b_2,b_3,b_4}} \prod_{i=1}^n \frac{1}{(z_i-z'_{\eta(i)})^2}.
\end{align}
Note that
\[
\prod_{i=1}^n \frac{1}{(z_i-z'_{\eta(i)})^2}=\prod_{\ell=1}^n \frac{1}{(\lambda_{\alpha_{i,1}}-\lambda_{\alpha_{i,2}})^2}
\]
where $(\alpha_{1,1},\alpha_{1,2}), \dots, (\alpha_{n,1},\alpha_{n,2})$ is a matching of $\lst{2n}$, i.e.,~an element of $\Pi_{2n}$. The term corresponding to a particular matching $\alpha\in \Pi_{2n}$ shows up in the sum in \eqref{eq:Cauchy_B4_2} exactly $2^n n! b_1! b_2! b_3! b_4!$ times. Indeed:
\begin{itemize}\itemsep=0pt
\item[--] for each pair $(\alpha_{i,1},\alpha_{i,2})$ in the matching we can choose which element will be a $z_i$ ($2^n$ choices), this identifies the index sets $\{\sigma(1),\dots,\sigma(n)\}$ and $\{\sigma(n+1),\dots,\sigma(2n)\}$,
\item[--] we can choose the values $\sigma(1), \dots, \sigma(n)$ by choosing one of the $n!$ permutations of the corresponding set, this will identify the lists $\cB_{\sigma,j}$, $1\le j\le 4$, and the \emph{sets} corresponding to the lists $\cB'_{\sigma,j}$, $1\le j\le 4$,
\item[--] we can choose the ordering of elements within the sets corresponding to $\cB'_{\sigma,j}$, $1\le j\le 4$, this can be done $b_1!b_2!b_3!b_4!$ ways, this completely identifies $\sigma$ and $\eta$.
\end{itemize}

Our calculations show that the sum in \eqref{eq:Cauchy_B4_2} is equal to
\[
\frac{n!}{b_1!b_2!b_3!b_4!} 2^n \sum_{\alpha\in \Pi_{2n}} \prod_{i=1}^n \frac{1}{(\lambda_{\alpha_{i,1}}-\lambda_{\alpha_{i,2}})^2}=\binom{n}{b_1,b_2,b_3,b_4} 2^n \Hf\left(\frac{\ind_{i\neq j}}{(\lambda_i-\lambda_j)^2}\right).
\]
This, together with \eqref{eq:C_Hf}, shows that the contribution of the diagonal terms is indeed equal to the expression on the right side of \eqref{eq:Cauchy_sum_1}.

To finish the proof of our lemma, we need to show that the off-diagonal terms in \eqref{eq:Cauchy_B4_1_2} cancel out. For this, it is enough to show that if $\eta\neq \tilde \eta$ are fixed permutations in $S_{b_1,b_2,b_3,b_4}$, then the following sum is equal to zero
\begin{align}\label{eq:Cauchy_B4_3}
 \sum_{\sigma\in S_{2n}} \prod_{i=1}^n \frac{1}{z_i-z'_{\eta(i)}}\cdot \frac{1}{z_i-z'_{\tilde \eta(i)}}.
\end{align}
To prove this statement, we first note that
the pair $\eta$, $\tilde \eta$ generates a permutation of $\lst{2n}$ with~${i\to n+\eta(i)}$, $n+i\to \tilde \eta^{-1}(i)$, $1\le i\le n$. This permutation has only even length cycles, and since $\eta\neq \tilde \eta$, it has at least one cycle with length $2\ell\ge 4$. Let one of these cycles~be
\[
i_1\to n+j_1\to i_2\to n+j_2\to \dots\to n+j_\ell\to i_1.
\]
The contribution of this cycle to the term corresponding to a particular $\sigma\in S_{2n}$ in \eqref{eq:Cauchy_B4_3} is
\begin{align}\label{eq:Cauchy_B4_4}
(-1)^k \prod_{k=1}^\ell \frac{1}{z_{i_k}-z'_{j_k}}\cdot \frac{1}{z'_{j_k}-z_{z_{i_{k+1}}}}.
\end{align}
(With $i_{n+1}=i_1$.)
The proof now follows along the line of the proof of \eqref{eq:off_diag_vanish_1}. Let
\[
\cF=\{i_1, \dots, i_\ell, n+j_1, \dots, n+j_\ell\}.
\]
We can evaluate the sum in \eqref{eq:Cauchy_B4_3} in two steps by
first fixing the values of $\sigma$ outside the index set $\cF$, and then assigning the values of
\[
\sigma(i_1), \dots, \sigma(i_\ell), \sigma(n+j_1), \dots, \sigma(n+j_\ell)
\]
out of the remaining $2\ell$ indices in all possible ways.
By Lemma~\ref{lem:cycle}, the sum of the terms \eqref{eq:Cauchy_B4_4} in this second step will be 0, which proves that the sum in \eqref{eq:Cauchy_B4_3} is also~0. This shows that the contribution of the off-diagonal terms in \eqref{eq:Cauchy_B4_1_2} vanishes, proving our lemma.
\end{proof}

\section{Asymptotics}
\label{sec:asymptotics}

We provide sketches of the proofs of Theorem~\ref{thm:limit_op} and Corollary \ref{cor:osc}, along with
the complementary Theorem~\ref{thm:limit_op1} and Corollary \ref{cor:osc1}.
For completeness, we start with the following statement regarding the limiting spectral measures of our block tridiagonal ensembles.

\begin{Theorem}\label{thm:eigcount}
Denote by $\mu_\gamma$ the semi-circle distribution with density \smash{$ \frac{1}{\pi \gamma} \sqrt{(4 \gamma - \lambda^2)_+}$}. For any integer $r$, any $s\ge 0$, the empirical spectral measure \smash{$\frac{1}{rn} \sum_{i=1}^{rn} \delta_{(r n)^{-1/2} \lambda_i}$} of the scaled eigenvalues of $\HH(r,s)$
converges weakly in law to $\mu_\gamma$ with $\gamma = \frac{r+s}{r}$. Likewise, for any integer $r$ and $s>0$ the suitably scaled empirical spectral measure of $\WW(r,s)$ converges weakly in law to
a scaled Marchenko--Pastur distribution as $n,m \rightarrow \infty$ with $m/n \rightarrow c \in [1, \infty)$.
\end{Theorem}

\begin{proof}
We follow the method of Trotter \cite{Trotter}, which
can be summarized thus.
Suppose that we have a sequence of symmetric $m\times m$ random matrices $A_m$ that are finitely banded, i.e., $[A_m]_{i,j}=0$ if $|i-j|>k$. Assume that it holds:
(i) $\frac1m E\|A_m-E A_m\|^2_2\to 0$ as $m \to \infty$, and (ii)
there exist $L^2[0,1]$ functions $h_0, \dots, h_k$ so that the $j$th diagonal in $EA_m$ for $0\le j\le k$ converges to $h_j$
after the natural embedding.
Then, the empirical spectral distribution of $A_m$ converges to the distribution of the random variable \smash{$\sum_{j=-k}^k h_{|j|}(U) \cos(2\pi j V)$} where $U$ and $V$
are independent uniform random variables on $[0,1]$.

Now consider $A_m = A_{rn}$ drawn from the law of \smash{$ \frac{1}
{\sqrt{rn}} \HH(r,s) $}. This is a symmetric band matrix, with $[A_{rn}]_{i,j}=0$ for $|i-j|> 2r$. We also have \smash{$\frac{1}{rn} E\|A_{rn}-EA_{rn}\|_2^2\to 0$} as there are~$\mathcal{O}(n)$ non-zero entries in $A_{rn}-EA_{rn}$, each of which has $\mathcal{O}(n^{-1})$ variance. By the properties of the scaled $\chi$
random variables appearing in the $r$th diagonal of $EA_{rn}$ we readily identify limiting functions $h_0, \dots, h_{2r}$ as zero except for
$
h_r(x)=\sqrt{\gamma}\sqrt{1-x}.
$
Hence, the given normalized empirical spectral distribution of $H_{n, \beta}{(r,s)}$ converges to the distribution of
\[Z =2 \sqrt{\gamma}\sqrt{1-U} \cos(2\pi r V) \stackrel{d}{=}
 2 \sqrt{\gamma}\sqrt{U'} \cos(2\pi V'),
\]
with new independent $[0,1]$-valued uniforms $U'$ and $V'$.
This though has the same distribution of the $x$ coordinate for a uniformly chosen point in $\{|z| \le 2 \sqrt{\gamma}\}$ which is exactly the advertised scaled semi-circle law.

Everything works much the same for a matrix drawn from the law $ \frac{1}{rn} \WW(r,s)$. This is again a $(2r)$-banded $rn \times rn$ matrix. As $n \rightarrow \infty$ with
$m/n \to c\in [1,\infty)$, Trotter's conditions $(i)$ and $(ii)$ hold with the only non-zero $h$-functions being $h_0(x)=\gamma(1+c-2x)$ and $h_r(x)=\gamma \sqrt{1-x}\sqrt{c-x}$.
This identifies the limiting empirical eigenvalue distribution in this case as that of then random variable
\[
W = \gamma(1+c-2U)+2 \gamma \sqrt{1-U}\sqrt{c-U} \cos(2\pi r V),
\]
with $U$ and $V$ as before. Recognizing the latter as a (scaled) Marchenko--Pastur law is not so immediate; the necessary details can be found in \cite[Proposition~4.2]{Ledoux}.
\end{proof}

\subsection[The H\_n,beta(r,s) soft edge]{The $\boldsymbol{{\tt{H}}_{n,\beta}(r,s)}$ soft edge}

Let $\mbf{T}_n$ be distributed as
${\tt{H}}_{n,\beta}(r,s)$ and consider the
centered and scaled matrix model
\begin{gather}
\mathbf{H}_n = \gamma^{-1/2} (rn)^{1/6} \bigl(2 \sqrt{(r+s)n} \mbf{I}_{rn} - \mathbf{T}_n \bigr)
 = 2 m_n^{2} \mbf{I}_{rn} - \sqrt{\frac{m_n}{{\gamma}}} \mbf{T}_n, \nonumber\\ m_n = (rn)^{1/3}.\label{eq:Hn}
\end{gather}
The introduced centering can be understood in light of Theorem~\ref{thm:eigcount}. We identify $m_n$ as a~continuum scale and view
 $v \in \mathbb{F}^{nr}$, on which $\mathbf{H}_n$ acts, as a vector-valued step function
\begin{equation}
\label{embed}
 (v_0, v_1, \dots) \in \ell_{nr}^2 \mapsto v(x) = v_{\lfloor m_n x \rfloor} \in L^2(\R_+, \mathbb{F}^r).
\end{equation}
Said another way, for integers $i\in [0, r-1] $ and $k \in [0,n]$,
$v_{i + k r}$ is identified with $v_i( k/m_{n})$
With this embedding in mind, write $
 \mathbf{H}_n = m_n^2 \Delta_n + m_n {\mathbf{V}}_n$.
Here $\Delta_n$ denotes discrete Laplacian in the underlying basis: it has diagonal blocks $2{I}_r$ and off-diagonal blocks
$- {I}_r$. The $\mbf{V}_n$ is then considered to define a matrix-valued potential.
The main technical result of
 \cite{Spike2} provides criteria for the convergence of the spectrum of $\mathbf{H}_n$ (in law)
in terms of the convergence of the integrated components of $\mathbf{V}_n$. A version of this, tailored to the particular setting considered here, is the following.

\begin{Proposition}
\label{Model_SoftEdgeConv}
Suppose that $\mbf{T}_n$ is an $rn\times rn$ Jacobi $r$-block matrix with blocks $\mbf{A}_{n,i}$, $\mbf{B}_{n,i}$. Let~$\mbf{H}_n$ and $m_n$ be defined as in \eqref{eq:Hn}.
Expressing $\mathbf{H}_n$ as $ m_n^2 \Delta_n + m_n {\mathbf{V}}_n$,
denote the running sums of the diagonal and upper diagonal block components of the potential $\mbf{V}_n$ by: for $0 \le x \le (rn)^{2/3}$,
\begin{align*}
\label{eq:potential_int}
 {\bf{Y}}_{n,1}(x) = \sum_{i=1}^{[m_n x]} \left( - \frac{1}{\sqrt{\gamma m_n}}
 \mathbf{A}_{n,i} \right),
 \qquad
 {\bf{Y}}_{n,2}(x) = \sum_{i=1}^{[m_n x]}
\left( m_n^2 \mbf{I}_r -
 \frac{1}{\sqrt{\gamma m_n}} \mathbf{B}_{n,i} \right).
\end{align*}
Now assume that
\begin{itemize}\itemsep=0pt
\item[$(i)$] Both $\{ {\bf{Y}}_{n,1}(x)\}_{x \ge 0}$ and
$\{ {\bf{Y}}_{n,2}(x)\}_{x \ge 0}$ are tight and
\begin{equation*}
\label{potential_convergence}
 {\bf{Y}}_{n,1}(x) +
 \bigl({\bf{Y}}_{n,2}(x) + {\bf{Y}}_{n,2}^{\dagger}(x) \bigr)\ \Rightarrow \
 \frac{1}{2} r x^2 \mbf{I}_r+ \sqrt{\frac{2}{\gamma}} B(x)
\end{equation*}
in law in the uniform-on-compacts topology;

\item[$(ii)$] For $j=1,2$, defining the increments via
\smash{$ {\bf{Y}}_{n,j}(x):=\frac{1}{\sqrt{\gamma m_n}} \sum_{i=1}^{[m_n x]}
 (\Delta {\bf{Y}})_{n,j,i}$}
 there is a~decomposition
\begin{equation*}
\label{potential_decomp}
 (\Delta {\bf{Y}})_{n, j,i}
 = \nn_{n,j,i} + (\Delta {\ww})_{n,j,i}
\end{equation*}
such that the following holds. The $\nn_{n,j, i} $ are diagonal and
$\nn_{n,j}(x) = \nn_{n,j,[m_nx]}$ satisfy
\begin{align}
 \kappa^{-1} x - \kappa \le \nn_{n,1}(x) + \nn_{n,2}(x) \le \kappa x + \kappa, \qquad
 \nn_{n,2}(x) \le 2 m_n, \label{compact1}
\end{align}
for some $\kappa \ge 1$. Also, with \smash{$\ww_{n,j}(x)=\sum_{i=1}^{[m_nx]} (\Delta {\ww})_{n,j,i}$}, there exists an $\epsilon > 0$ so that
\begin{align}
 \| \ww_{n,1}(x) - \ww_{n,1}(y)\|^2 + \| \ww_{n,2}(x) - \ww_{n,2}(y)\|^2 & \le \kappa_n x^{1-\epsilon} + \kappa_n, \label{compact3}
\end{align}
for all $|x-y| \le 1$. Here
$\|\cdot\|$ denotes the spectral norm, and $\kappa_n \ge 1$ a sequence of tight random constants.
\end{itemize}

Then, with $\mathcal{H}_{\beta, \gamma}$ defined in \eqref{eq:H_op},
we have
that any finite collection of ordered eigenvalues of~$\mathbf{H}_n$, along with their associated eigenfunctions as elements in $L^2 (\mathbb{F}^r)$, converge jointly in law to the corresponding eigenvalues/eigenfunctions of $ \mathcal{H}_{\beta, \gamma}$.
\end{Proposition}

Condition (i) simply identifies the correct limit potential. Condition (ii) provides an almost sure lower bound on $\langle v, \mathbf{H}_n v \rangle$ independent of $n$ which is essential for extracting eigenvalue limits. The bound \eqref{compact3} controls the random oscillations of the potential by the growth of $\nn_{n,j}$, which by \eqref{compact1} is well-controlled in terms of its natural limit.

Both conditions are readily checked in our case. The independence of the entries of
${\bf{Y}}_{n,1}$ and~${\bf{Y}}_{n,2}$ allows one to establish the functional limit theorem in (i) component-wise. For (ii) one can take \smash{$\nn_{n,0} =0$} and \smash{$\nn_{n,2}$} the vector of expectations of centered $\chi$ variables appearing on the diagonal of ${\bf{Y}}_{n,2}$.
In both (i) and (ii) that the sequence $\chi_p -\sqrt{p}$ converges in
law to a $N\bigl(0,\frac{1}{2}\bigr)$ random variable as $p \uparrow \infty$ and satisfies a uniform in $p$ subgaussian bound plays a~fundamental role.
With the matrix models considered in \cite{Spike2} so similar to $\mathbf{H}_n$, the necessary details are effectively identical to what has already been done there.

The proof of Corollary \ref{cor:osc} is also similar to that of the
corresponding statement in \cite{Spike2}, though here is a sketch of the main idea. Fix $\lambda\in \R$ and consider the system
\begin{align}
{\rm d}F(x) = F'(x) {\rm d}x, \qquad
 {\rm d}F'(x) = (\lambda - r x) F(x) + \sqrt{\frac{2}{\gamma}} F (x) {\rm d}B(x),\label{Eig:system}
\end{align}
with initial condition $F$ satisfying $(F(0), F'(0)) = (0, \mbf{I}_r)$. Here $F$ is an $r\times r$ matrix valued function.
It can be shown that the number of eigenvalues of
the Dirichlet problem for $\mathcal{H}_{\beta, \gamma}$ less than $\lambda$ coincides with the number of zeros of $\det F$ on $\R_+$

Next define $P(x) = F'(x) F(x)^{-1} $, the matrix Riccati substitution. This satisfies
\begin{equation}
\label{matrix_ricatti}
{\rm d} P(x) = \bigl((\lambda - r x)\mbf{I}_r - P^2(x)\bigr){\rm d}x + \sqrt{\frac{2}{\gamma}} {\rm d}B(x).
\end{equation}
One may then verify that points $x'$ where $\det F$ vanishes correspond to $P(x)$ possessing an eigenvalue $p(x)$ that explodes to $-\infty$ as $x \rightarrow x'$. This is exactly the content of Theorem~\ref{cor:osc}: the stochastic differential equation \eqref{mult_sde}, as can be derived from \eqref{matrix_ricatti} by an application of It\^{o}'s Lemma, describes the evolution of the eigenvalues
$(p_1(x), \dots, p_r(x))$ of $P(x)$, continued through explosion times to $-\infty$. (Again, see \cite{Spike2} for additional details.)

\subsection[Hard edge for W\_n, n+a,beta(r,s)]{Hard edge for $\boldsymbol{{\tt{W}}_{n, n+a,\beta}(r,s)}$}

To give a precise definition of
 $\mathcal{G}_{\beta, \gamma}$ from \eqref{matrixgenerator}, we first have to define $\mbf{Z}_x$. We introduce a new type of $r \times r$ matrix
Brownian motion $x \mapsto B_x$ in which all entries are independent, the off-diagonal entries are standard ${\mathbb{F}}$-Brownian motions and the diagonal entries are real Brownian motions with common diffusion coefficient $\frac{1}{\beta}$. Then the coefficient matrix $\mbf{Z}_x$
is given by
\begin{equation}
\label{WandA}
 \mbf{Z}_x = \mbf{Y}_x \mbf{Y}_x^{\dagger}, \qquad
 \mbf{Y}_x^{-1} {\rm d} \mbf{Y}_x = \frac{1}{\sqrt{\gamma}}{\rm d}B_x +
 \left( \frac{a}{2 \gamma} - \frac{1}{2\beta \gamma} \right) \mbf{I}_r {\rm d}x,
\end{equation}
in which $\mbf{Y}_0 = 0.$ Notice that in the $r=1$ and $\gamma =1$ setting in which
the Stochastic Bessel Operator was first introduced, \smash{$Z_x =
{\rm e}^{ \frac{2}{\sqrt{\beta}} b(x) + a x}$} with a standard one-dimensional Brownian motion~$b(x)$.

While \eqref{matrixgenerator} is a nice format in which to package the limiting operator $\mathcal{G}_{\beta, \gamma}$, we actually identify this operator via its inverse. An exercise shows that, specifying a Dirichlet condition at the origin, $\mathcal{G}_{\beta, \gamma} = \mathcal{L}_{\beta, \gamma} \mathcal{L}_{\beta, \gamma}^{\dagger} $,
where
\begin{equation*}
\label{kernel_ops}
 \mathcal{L}_{\beta, \gamma}^{-1} f (x) = \int_x^\infty {\ell}(x,y) f(y) {\rm d}y, \qquad {\ell}(x,y) = {\rm e}^{-rx/2} \mbf{Y}_x^{-1} \mbf{Y}_y.
\end{equation*}
The results of \cite{RR2} imply that $\mathcal{L}_{\beta, \gamma}^{-1}$ is Hilbert--Schmidt for any $\gamma >0$ and $a>-1$, and the main convergence result there can be summarized, in the spirit of Proposition~\ref{Model_SoftEdgeConv}, as follows.

\begin{Proposition}\label{prop:hard_lim}
Given a block bidiagonal matrix $\mbf{L}_n$ with independent diagonal and upper diagonal $r \times r$ blocks $\mbf{D}_k=\mbf{D}_{n,k}$ and $\mbf{O}_k=\mbf{O}_{n,k}$, embed $ \mbf{ L}_n^{-1}$ into $L^2([0,1], \mathbb{F}^r)$
in the manner of~\eqref{embed} with now $m_n = n$. In particular,
define the piecewise step kernel
\begin{equation*}
\label{explicit_inv}
 \ell_{n}(x,y) = \sqrt{ \frac{ {n} \gamma}{r}}
 {\mbf{D}}_{\lfloor nx \rfloor}^{-1}
 {\mbf{X}}_n(x)^{-1} \mbf{X}_n(y) \mathbf{1}_{0 \le x < y \le 1},
 \qquad
 {\mbf{X}}_n(x) = \prod_{k= 1}^{ \lfloor nx \rfloor}
 {\mbf{O}}_k {\mbf{D}}_{k+1}^{-1}.
\end{equation*}
Assume that
\begin{itemize}\itemsep=0pt
\item[$(i)$] As $n \rightarrow \infty$, in the uniform-on-compacts topology
\[
 \left( \sqrt{\frac{{n} \gamma}{r}} \mbf{D}_{\lfloor nx \rfloor}^{-1}, \mbf{X}_n(x) \right) \Rightarrow
 \left( \frac{1}{r \sqrt{ 1-x}} \mbf{I}_r, \mbf{X}_x \right) \qquad
\text{for} \ x\in [0,1),
\]
where $\mbf{X}_x$ satisfies the matrix stochastic differential equation
\[
 \mbf{X}_x^{-1} {\rm d}\mbf{X}_x = \frac{1}{\sqrt{r \gamma (1-x)}} {\rm d} B_x - \frac{a -\frac{1}{ \beta}}{2 r\gamma (1- x)} {I}_r{\rm d}x, \qquad \mbf{X}_0=0.
\]

\item[$(ii)$] There is the bound
\[
 \int_0^1 \int_0^y | \ell_n(x,y) |^2 {\rm d}x {\rm d}y \le \kappa_n
\]
with $\kappa_n$ a sequence of tight random constants.
\end{itemize}
Then $\operatorname{spec}\bigl( \frac{rn}{ \gamma} \mbf{L}_n \mbf{L}_n^\dagger\bigr) \rightarrow \operatorname{spec}(\mathcal{G}_{\beta, \gamma})$ in the manner described in Proposition {\rm\ref{Model_SoftEdgeConv}}.
\end{Proposition}

The point is that the integral kernel $\ell_n(x,y)$ is an exact representation of \smash{$\bigl(\sqrt{rn / \gamma} \mbf{L}_n\bigr)^{-1}$}.
The convergence in (i) identifies the pointwise limit of $\ell_n$ while (ii) implies (subsequential) convergence of the corresponding integral operator in Hilbert--Schmidt norm.
 Combined, we have that \smash{$ \bigl(\frac{rn}{ \gamma} \mbf{L}_n \mbf{L}_n^{\dagger}\bigr)^{-1}$} converges in norm in the same subsequential
sense; convergence of the finite parts of the spectrum follows.

Checking the conditions of Proposition~\ref{prop:hard_lim} for ${\tt{W}}_{n, n+a,\beta}(r,s)$ follows the arguments of \cite{RR2}.
For (i), by an explicit expansion of the inverse of the diagonal blocks, the increments of
$\mbf{X}_n$ have the form
\[
 {\mbf{O}}_k {\mbf{D}}_{k+1}^{-1} = \mbf{I}_r + \frac{1}{\sqrt{ (r+s)(n-k)}} \mbf{G}_k + \frac{-a+ \frac{1}{\beta} }{2(r+s) (n-k)} \mbf{I}_r + \boldsymbol{\varepsilon}_k.
\]
Here the $\mbf{G}_k$ are independent and have independent entries which are $\FF$-normals off diagonal and centered $\chi$ variables on the diagonal. As such, the $\mbf{G}_k$ form the increments of the limiting matrix Brownian motion $B_x$. The
$\boldsymbol{\varepsilon}_k$ matrices are error terms for which
one can derive the estimate $E \|\boldsymbol{\eps}_k \|^p
= \mathcal{O}\bigl( (n-k)^{-\frac{3}{2} p}\bigr)$.
The proof of (ii) is far more technical, as it requires sharp control of the paths of $x \mapsto \mbf{X}_n(x)$ in a vicinity of $x=1$.

Note that while the embedding of the matrix $\mbf{L}_n^{-1}$ naturally takes place as an operator on $[0,1]$, the advertised limit $\mathcal{L}^{-1}_{\beta, \gamma}$ lives on $[0, \infty)$.
This is just more convenient for eventual comparison to the soft edge operator
$\mathcal{H}_{\beta, \gamma}$, and the kernel $\ell(x,y)$ is related to the limit of $\ell_n$, which is constructed from the limit objects defined in (i) above, by the simple change of variables $x \mapsto 1 - {\rm e}^{-rx} =\varphi(x)$. As one can readily check, $ \mbf{X}_{\varphi(x)} = \mbf{Y}_x$ and
\smash{$ \frac{1}{r\sqrt{1-\varphi(x)}} {\rm d}\varphi(x) = {\rm e}^{-rx/2} {\rm d}x$}. This establishes Theorem~\ref{thm:limit_op1}

From here, the proof of Corollary \ref{cor:osc1} amounts to writing out the eigenvalue problem for $\mathcal{G}_{\beta, \gamma}$ as a system, then invoking the matrix Riccati correspondence, in analogy with equations \eqref{Eig:system} and \eqref{matrix_ricatti}. \cite[Section 4]{RR2} provides the details.

We conclude by recording the hard edge versions of Corollary \ref{cor:betalimit} and Conjecture \ref{con:betalimit}. Denote by
$\operatorname{Bessel}_{\beta, a}$ the random point process defined by the
$r=1$ and $\gamma =1$ case of Theorem~\ref{thm:limit_op1} (and Corollary \ref{cor:osc1}). To be totally concrete, the points of $\operatorname{Bessel}_{\beta, a}$
are the Dirichlet eigenvalues of the one-dimensional operator
\begin{align*}
 -{\rm e}^{(a+1)x +\frac{2}{\sqrt{\beta}} b(x)} \frac{\rm d}{{\rm d}x} {\rm e}^{-ax -\frac{2}{\sqrt{\beta}} b(x)} \frac{\rm d}{{\rm d}x} \cdot
\end{align*}
 acting on the positive half-line. This is well defined
for all $\beta > 0$ and $a> -1$.

\begin{Corollary}
Lut us consider the eigenvalue point processes for our solvable instances of
${\tt{W}}_{n, n+a,\beta}(r, s)$. These have explicit joint densities proportional to
 \begin{align}
 \label{density3}
 |\Delta({\lambda})|^{\beta}
 \left( \sum_{(\mathcal{A}_1,\dots,\mathcal{A}_r)\in \cP_{r,n}} \prod_{j=1}^r \Delta(\cA_j)^2 \right) \prod_{i=1}^{rn} \lambda_i^{\frac{\beta}{2}( (r+s)a+1)-1} {\rm e}^{-\frac{\beta}{2} \lambda_i}
 \ind_{\R_n^+}
 \end{align}
for $r \ge 2$ and $\beta s=2$, and to
 \begin{align}
 \label{density4}
 \Delta({\lambda})^{\beta+\frac{\beta s}{2}} \left|\Pf \left(\frac{{\bf{1}}_{i \neq j}}{\lambda_i -\lambda_j} \right)\right|^{\frac{\beta s}{2}}
 \prod_{i=1}^{rn} \lambda_i^{\frac{\beta}{2}( (r+s) a+1)-1} {\rm e}^{-\frac{\beta}{2} \lambda_i} \ind_{\R_n^+}
 \end{align}
for $r=2$ and $\beta s =2$ or $4$. When $r=2$, $\beta s =2 $ and $\beta=1$, the scaling limit of the minimal points is given by ${\operatorname{Bessel}}_{2, a/2}$. When $r=2$, $\beta s = 4$, and $\beta = 2$, the scaling limit is $\operatorname{Bessel}_{4, a/2}$.
\end{Corollary}

\begin{Conjecture}
More generally, the minimal points under \eqref{density3} have scaling limit given by \smash{${\operatorname{Bessel}}_{\beta +\frac{2}{r}, a/(1+\frac{2}{\beta r})}$} for
$r \ge 2$ and $\beta =1$ or $2$. For the minimal points under \eqref{density4} with~${\beta s =4}$ and $\beta =1$, the scaling limit is instead~\smash{${\operatorname{Bessel}}_{3, a/3}$}.
\end{Conjecture}

\subsection*{Acknowledgements}

The authors thank Philippe Di Francesco for pointing out reference \cite{DSZ}. Thanks as well to the anonymous referees for a number of observations and useful suggestions. B.V.~was partially supported by the University of Wisconsin – Madison Office of the Vice Chancellor for Research and Graduate Education with funding from the Wisconsin Alumni Research Foundation and by the National Science Foundation award DMS-2246435.

\pdfbookmark[1]{References}{ref}
\LastPageEnding

\end{document}